%&amstex          
\input amstex\documentstyle{amsppt}  
\pagewidth{12.5cm}\pageheight{19cm}\magnification\magstep1
\topmatter
\title Character sheaves on disconnected groups, IX\endtitle
\author G. Lusztig\endauthor
\address{Department of Mathematics, M.I.T., Cambridge, MA 02139}\endaddress
\thanks{Supported in part by the National Science Foundation.}\endthanks
\endtopmatter   
\document

\define\ufs{\un{\fs}}
\define\Bpq{\Bumpeq}

\define\tir{\ti{\rho}}

\define\ucf{\un{\cf}}

\define\du{\dot u}
\define\dw{\dot w}

\define\uD{\un D}

\define\mpb{\medpagebreak}
\define\co{\eta}
\define\Up{\Upsilon}

\define\bd{\bar d}
\define\bh{\bar h}

\define\bL{\bar L}

\define\btC{\bar{\tC}}
\define\bPP{\bar{\PP}}

\define\hD{\hat D}
\define\hH{\hat H}
\define\hZ{\hat Z}

\define\dsv{\dashv}

\define\po{\text{\rm pos}}

\define\si{\sim}

\define\sqc{\sqcup}

\define\tcl{\ti\cl}
\define\tca{\ti\ca}

\define\tiz{\ti\z}

\define\bj{\bar j}

\define\bK{\bar K}
\define\bC{\bar C}

\define\bX{\bar X}
\define\bZ{\bar Z}

\define\bsi{\bar\s}
\define\lb{\linebreak}

\define\bin{\binom}
\define\op{\oplus}

\define\part{\partial}
\define\em{\emptyset}

\define\n{\notin}
\define\iy{\infty}
\define\m{\mapsto}
\define\do{\dots}

\define\sub{\subset}    
\define\bxt{\boxtimes}
\define\T{\times}
\define\ti{\tilde}
\define\nl{\newline}
\redefine\i{^{-1}}

\define\un{\underline}
\define\ov{\overline}
\define\ot{\otimes}
\define\bbq{\bar{\QQ}_l}

\define\Ad{\text{\rm Ad}}
\define\Hom{\text{\rm Hom}}

\define\IND{\text{\rm IND}}
\define\ind{\text{\rm ind}}

\define\tr{\text{\rm tr}}

\define\a{\alpha}
\redefine\b{\beta}
\redefine\c{\chi}
\define\g{\gamma}
\redefine\d{\delta}
\define\e{\epsilon}
\define\et{\eta}
\define\io{\iota}
\redefine\o{\omega}
\define\p{\pi}
\define\ph{\phi}
\define\ps{\psi}
\define\r{\rho}
\define\s{\sigma}
\redefine\t{\tau}

\define\k{\kappa}
\redefine\l{\lambda}
\define\z{\zeta}
\define\x{\xi}

\define\vt{\vartheta}

\redefine\D{\Delta}

\redefine\L{\Lambda}
\define\Ph{\Phi}

\define\boc{\bold c}
\define\dd{\bold d}

\define\kk{\bold k}

\define\rr{\bold r}

\define\EE{\bold E}

\define\II{\bold I}

\define\LL{\bold L}

\define\NN{\bold N}

\define\PP{\bold P}
\define\QQ{\bold Q}
\define\RR{\bold R}

\define\TT{\bold T}

\define\WW{\bold W}
\define\ZZ{\bold Z}
\define\XX{\bold X}

\define\ca{\Cal A}
\define\cb{\Cal B}

\define\cd{\Cal D}
\define\ce{\Cal E}
\define\cf{\Cal F}

\define\ch{\Cal H}

\define\ck{\Cal K}
\define\cl{\Cal L}
\define\cm{\Cal M}

\define\cp{\Cal P}

\define\cv{\Cal V}

\define\fe{\frak e}
\define\ff{\frak f}

\define\fh{\frak h}

\define\fj{\frak j}

\define\fs{\frak s}

\define\fK{\frak K}

\define\fU{\frak U}

\define\tc{\ti c}

\define\te{\ti e}
\define\tf{\ti f}

\define\tih{\ti h}
\define\tj{\ti j}

\define\tio{\ti\o}

\define\tit{\ti t}

\define\tA{\ti A}
\define\tB{\ti B}
\define\tC{\ti C}

\define\tK{\ti K}

\define\tT{\ti T}

\define\sh{\sharp}

\define\bc{\bar c}

\define\hC{\hat C}

\define\che{\check}
\define\cha{\che{\a}}

\define\bat{\bar\t}
\define\ucl{\un\cl}
\define\ucm{\un\cm}
\define\dcl{\dot{\cl}}
\define\udcl{\un{\dot{\cl}}}
\define\udcm{\un{\dot{\cm}}}

\define\ndsv{\not\dsv}
\define\prq{\preceq}
\define\BBD{BBD}
\define\DE{D}
\define\GR{Gr}
\define\CS{L3}
\define\AD{L9}
\define\PCS{L10}
\define\CRG{L14}
\define\MV{MV}
\define\GI{Gi}
\head Introduction\endhead
Throughout this paper, $G$ denotes a fixed, not necessarily connected, reductive algebraic
group over an algebraically closed field $\kk$. This paper is a part of a series \cite{\AD}
which attempts to develop a theory of character sheaves on $G$. 

One of the main constructions in \cite{\CS} (going back to \cite{\CRG}) was a procedure 
which to any character sheaf on $G^0$ associates a certain two-sided cell in an (extended)
Coxeter group. A variant of this construction (restricted to "unipotent" character sheaves)
was later given by Grojnowski \cite{\GR}. Here we give a construction which generalizes 
that in \cite{\CS} (and takes into account the approach in \cite{\GR}) which to any 
(parabolic) character sheaf on $Z_{J,D}$ associates a certain type of two-sided cell.

The paper is organized as follows. In Section 40 we study certain equivariant sheaves on  
$G^0/U^*\T G^0/U^*$ (where $U^*$ is the unipotent radical of a Borel in $G^0$) under the
convolution operation. Some results in this section are implicit in \cite{\CRG, Ch.1}. In 
Section 41 we study the character sheaves on $Z_{\em,D}$ (where $D$ is a connected
component of $G$) by connecting them with sheaves on $G^0/U^*\T G^0/U^*$. We use this study
to attach a two-sided cell to any character sheaf on $Z_{J,D}$. (See 41.4.) In Section 42 
we study the interaction between the duality operation $\dd$ (see 38.10, 38.11) and the 
functor $\ff_{\em,\II}$ (see 36.4). The main result in this section is Proposition 42.9 
which contains \cite{\CS, III, Cor.15.8(b)} as a special case (with $G=G^0,v=1$).

{\it Notation} We fix a $1$-dimensional $\bbq$-vector space $V$ with a given isomorphism 
$V^{\ot2}@>\si>>\bbq(1)$ (Tate twist of $\bbq$). For $n\in\NN$ we set 
$\bbq(n/2)=V^{\ot n}$. For $n\in\ZZ,n<0$ let $\bbq(n/2)$ be the dual space of $\bbq(-n/2)$.
If $X$ is an algebraic variety and $A\in\cd(X),n\in\ZZ$ we write $A[[n/2]]$ instead of 
$A[n](n/2)$. (When $n$ is even this agrees with the notation in \cite{\AD, II, p.73}.)

\head Contents \endhead
40. Sheaves on $G^0/U^*\T G^0/U^*$.

41. Character sheaves and two-sided cells.

42. Duality and the functor $\ff_{\em,\II}$.

\head 40. Sheaves on $G^0/U^*\T G^0/U^*$ \endhead
\subhead 40.1\endsubhead
Let $\ca=\ZZ[v,v\i]$. Let $\hH$ (resp. $H$) be the $\ca$-module consisting of all formal 
(resp. finite) linear combinations $\sum_{w\in\WW,\l\in\ufs}a_{w,\l}\tT_w1_\l$ with
$a_{w,\l}\in\ca$. Note that $H$ is naturally an $\ca$-submodule of $\hH$ with $\ca$-basis
$\{\tT_w1_\l;w\in\WW,\l\in\ufs\}$. For any $n\in\NN^*_\kk$, the $\ca$-submodule of $H$ 
spanned by $\{\tT_w1_\l;w\in\WW,\l\in\ufs_n\}$ may be naturally identified with $H_n$ (see
31.2(a)). There is a unique $\ca$-algebra structure on $\hH$ in which the product of two 
elements

$h=\sum_{w\in\WW,\l\in\ufs}a_{w,\l}\tT_w1_\l$,
$h'=\sum_{w'\in\WW,\l'\in\ufs}a'_{w',\l'}\tT_{w'}1_{\l'}$
\nl
as above is $hh'=\sum_{y\in\WW,\nu\in\ufs}b_{y',\nu}\tT_y1_\nu$ where for any $\nu\in\ufs$,

$\sum_{w,w'\in\WW}a_{w,w'{}\i\nu}a'_{w',\nu}\tT_w\tT_{w'}1_\nu=
\sum_{y\in\WW}b_{y,\nu}\tT_y1_\nu$
\nl
is computed in the algebra structure of $H_n$ for any $n$ such that $\nu\in\ufs_n$. Thus 
$\hH$ becomes an associative algebra with $1$; $H$ is a subalgebra (without $1$) and, for
$n\in\NN^*_\kk$, $H_n$ is a subalgebra (with a different $1$) with the algebra structure as
in 31.2.

Now in the definition of $\hH$ given above, although $\tT_w1_\l$ is defined, the elements 
$\tT_w,1_\l$ are not defined separately (as was the case in $H_n$). To remedy this we set 
$\tT_w=\sum_{\l\in\ufs}\tT_w1_\l\in\hH$ (for $w\in\WW$) and $1_\l=\tT_11_\l\in H$ (for
$\l\in\ufs$). Then $\tT_w1_\l$ is the product of $\tT_w,1_\l$ in the algebra $\hH$. Note
that $\tT_1$ is the unit element of $\hH$ and the following equalities hold in $\hH$:

$1_\l1_\l=1_\l$ for $\l\in\ufs,1_\l1_{\l'}=0$ for $\l\ne\l'$ in $\ufs$;

$\tT_w\tT_{w'}=\tT_{ww'}$ for $w,w'\in\WW$ such that $l(ww')=l(w)+l(w')$;

$\tT_w1_\l=1_{w\l}\tT_w$ for $w\in\WW,\l\in\ufs$;

$\tT_s^2=\tT_1+(v-v\i)\sum_{\l\in\ufs;s\in\WW_\l}\tT_s1_\l$ for $s\in\II$.
\nl
By a standard argument we see that 

(a) $H$ is exactly the $\ca$-algebra defined by the generators $\tT_w1_l$ ($w\in\WW$,
$\l\in\ufs$) and the relations:

$(\tT_w1_\l)(\tT_{w'}1_{\l'})=0$ if $w,w'\in\WW,\l,\l'\in\ufs,w'\l'\ne\l$,

$(\tT_w1_{w'\l'}(\tT_{w'}1_{\l'})=\tT_{ww'}1_{\l'}$ if $w,w'\in\WW,\l,\l'\in\ufs$,
$l(ww')=l(w)+l(w')$,

$(\tT_s1_{s\l'})(\tT_s1_{\l'})=\tT_11_{\l'}+(v-v\i)c\tT_s1_{\l'}$ if $s\in\II$,
$\l'\in\ufs$ where $c=1$ for $s\in\WW_{\l'}$ and $c=0$ for $s\n\WW_{\l'}$.

\subhead 40.2\endsubhead
Let $R,R^+$ be as in 28.3. The following result is well known:

(a) {\it If $w\in\WW,\a\in R^+$ and $s_\a$ is as in 28.3 then we have $l(ws_\a)>l(w)$ if 
and only if $w(\a)\in R^+$.}

Let $\l\in\ufs$. Let $R_\l,R^+_\l,\WW_\l,H_\l$ be as in 34.2. We write $\O_\l$ instead of 
$\O_\l^D$ (as in 34.4 with $D=G^0$). We show:

(b) {\it If $w\in\WW$ then $w\WW_\l$ contains a unique element $w_1$ of minimal length; it
is characterized by the property $w_1(R^+_\l)\sub R^+$.}
\nl
Let $w_1$ be an element of minimal length in $w\WW_\l$. Let $\a\in R_\l^+$. Then 
$l(w_1s_\a)\ge l(w_1)$. Since $l(w_1s_\a)=l(w_1)+1\mod 2$ we see that $l(w_1s_\a)>l(w_1)$.
By (a) we have $w_1(\a)\in R^+$. Thus, $w_1(R_\l^+)\sub R^+$. Now let $u\in\WW_\l-\{1\}$. 
We pick $\a\in R^+_\l$ such that $u(\a)\i\in R^+_\l$; then $w_1u(\a)\i\in R^+$. If $w_1u$ 
has minimal length in $w\WW_\l$ then by an earlier part of the argument applied to $w_1u$ 
instead of $w_1$ we have $w_1u(\a)\in R^+$, a contradiction. We see that $w_1$ is the 
unique element of minimal length in $w\WW_\l$. It remains to show that if $u\in\WW_\l$ 
satisfies $w_1u(R^+_\l)\sub R^+$ then $u=1$. If $u\ne 1$ then by an earlier part of the 
argument we have $w_1u(\a)\i\in R^+$ for some $\a\in R^+_\l$, a contradiction. This proves
(b).

We show:

(c) {\it If $s\in\II$ and $w\in\WW$ has minimal length in $w\WW_\l$ then either (i) $sw$
has minimal length in $sw\WW_\l$ or (ii) $w\i sw\in\WW_\l$.}
\nl
There is a unique $\b\in R^+$ such that $s(\b)\i\in R^+$. Assume that (i) does not hold. By
(b) there exists $\a\in R^+_\l$ such that $sw(\a)\i\in R^+$; moreover, $w(\a)\in R^+$. 
Hence $w(\a)=\b$. We have $w\i(\b)=\a\in R_\l$ hence $w\i sw\in\WW_\l$ and (ii) holds. This
proves (c).

For $z\in\WW_\l$ let $\tT^\l_z,c^\l_z\in H_\l$ be as in 34.2 . Then 
$c^\l_z=\sum_{z'\in\WW_\l}p^\l_{z',z}\tT^\l_{z'}$ where $p^\l_{z',z}\in\ZZ[v\i]$ are
uniquely defined.

For any $w\in\WW$, $\l\in\ufs$ there is a unique element element of $H$ which is equal to 
$c_{w,\l}\in H_n$ (see 34.4) for any $n$ such that $\l\in\ufs_n$; we denote this element
again by $c_{w,\l}$.  We have 

$c_{w,\l}=\sum_{w'\in\WW}\p_{w',w,\l}\tT_{w'}1_\l$
\nl
where $\p_{w',w,\l}\in\ZZ[v\i]$ are uniquely defined. Note that 

{\it $\{c_{w,\l};w\in\WW,\l\in\ufs\}$ is an $\ca$-basis of $H$.}

We show:

(d) {\it Let $w,w'\in\WW$. We write $w=w_1z,w'=w'_1z'$ where $w_1$ has minimal length in
$w\WW_\l$, $w'_1$ has minimal length in $w'\WW_\l$ and $z,z'\in\WW_\l$. If $w_1\ne w'_1$ 
then $\p_{w',w,\l}=0$. If $w_1=w'_1$ then $\p_{w',w,\l}=p^\l_{z',z}$.}
\nl
From the definitions we see that if $w\l\ne w'\l$ then $\p_{w',w,\l}=0$. Thus we may assume
that $w\l=w'\l$. We choose a sequence $s_1,s_2,\do,s_r$ in $\II$ such that
$w\l=w'\l=s_rs_{r-1}\do s_1\l\ne s_{r-1}\do s_1\l\ne\do\ne s_1\l\ne\l$.

We show that for $k\in[0,r]$, $s_ks_{k-1}\do s_1$ has minimal length in 
$s_ks_{k-1}\do s_1\WW_\l$. We argue by induction. For $k=0$ the result is obvious. Assume 
now that $k\in[1,r]$. Since $s_{k-1}\do s_1$ has minimal length in $s_{k-1}\do s_1\WW_\l$ 
and $s_ks_{k-1}\do s_1\l\ne s_{k-1}\do s_1\l$ we see from (c) that $s_ks_{k-1}\do s_1$ has
minimal length in $s_ks_{k-1}\do s_1\WW_\l$ as required.

In particular, $s_rs_{r-1}\do s_1$ has minimal length in $s_rs_{r-1}\do s_1\WW_\l$. Since 
$w\l=s_rs_{r-1}\do s_1\l$ we have $w=s_rs_{r-1}\do s_1h_1h_2$ where $h_1\in\O_\l$,
$h_2\in\WW_\l$. Then both $w_1$ and $s_rs_{r-1}\do s_1h_1$ have minimal length in 
$s_rs_{r-1}\do s_1h_1\WW_\l=w\WW_\l=w_1\WW_\l$; using (b) we deduce that 
$s_rs_{r-1}\do s_1h_1=w_1$. Hence $s_1\do s_rw=s_1\do s_rw_1z=h_1z$. Similarly,
$s_1\do s_rw'=h'_1z'$ where $h'_1\in\O_\l$.

From the results in 34.7-34.10 we see that 
$\p_{w',w,\l}=p^\l_{s_1\do s_rw',s_1\do s_rw}=p^\l_{h'_1z',h_1z}$. Using $h_1,h'_1\in\O_\l$
and the definitions (34.2) we see that $p^\l_{h'_1z',h_1z}=0$ if $h_1\ne h'_1$ and 
$p^\l_{h'_1z',h_1z}=p^\l_{z',z}$ if $h_1=h'_1$.

It remains to show that we have $w_1=w'_1$ if and only if $h_1=h'_1$. We have 
$s_rs_{r-1}\do s_1=h_1\i w_1$ and similarly $s_rs_{r-1}\do s_1=(h'_1)\i w'_1$. Hence 
$h_1\i w_1=(h'_1)\i w'_1$. We see that $w_1=w'_1$ if and only if $h_1=h'_1$. This proves
(d).

\mpb

For $w'\le w$ in $\WW$, $\l\in\ufs$ and $i\in\ZZ$ we define $N_{i,w',w,\l}\in\ZZ$ by
 
(e) $\p_{w',w,\l}=v^{l(w')-l(w)}\sum_{i\in\ZZ}N_{i,w',w,\l}v^i$, that is,

$p^\l_{z',z}=v^{l(w')-l(w)}\sum_{i\in\ZZ}N_{i,w',w,\l}v^i$ if $w'\WW_\l=w\WW_\l$ and
$z,z'$ are as in (d), 

$N_{i,w',w,\l}=0$ if $w'\WW_\l\ne w\WW_\l$.
\nl
Note that $N_{i,w',w,\l}$ is $0$ unless $i$ is even.

\subhead 40.3\endsubhead
Let $B^*\in\cb$. Let $U^*=U_{B^*}$ and let $T$ be a maximal torus of $B^*$. Let 
$\rr=\dim\TT$. Let $W_T=N_{G^0}T/T$. We identify $T=\TT,W_T=\WW$ as in 28.5. For any 
$w\in\WW$ we denote by $\dw$ a representative of $w$ in $N_{G^0}T$.

Let $C=G^0/U^*\T G^0/U^*$. We have a partition $C=\cup_{w\in\WW}C_w$ where

$C_w=\{(hU^*,h'U^*)\in C;h\i h'\in B^*\dw B^*\}$.
\nl
For $w\in\WW$ let $d_w=\dim C_w$ and let

$\bC_w=\{(hU^*,h'U^*)\in C;h\i h'\in\ov{B^*\dw B^*}\}$
\nl
(closure in $G^0$). Now $\bC_w$ is an irreducible variety and we have a partition 
$\bC_w=\cup_{w';w'\le w}C_{w'}$ with $C_w$ smooth, open dense in $\bC_w$.

Define $\g_{\dw}:B^*\dw B^*@>>>T$ by $\g_{\dw}(g)=t$ where $g\in U^*\dw tU^*$ with 
$t\in T$. Define $\ps:C_w@>>>T$ by $\ps(hU^*,h'U^*)=\g_{\dw}(h\i h')$. 

For $\cl\in\fs$ we set $\cl_w=\ps^*\cl$, a local system on $C_w$. (Using 28.1(c) we see 
that the isomorphism class of $\ps^*\cl$ is independent of the choice of $\dw$.) Let 
$\cl_w^\sh=IC(\bC_w,\cl_w)\in\cd(\bC_w)$. 

\subhead 40.4\endsubhead
For $w\in\WW,\cl\in\fs$ let $\ucl_w=j_{w!}\cl_w$, $\ucl^\sh_w=\bj_{w!}\cl^\sh_w$ where 
$j_w:C_w@>>>C$, $\bj_w:\bC_w@>>>C$ are the inclusions. Let $\hC$ be the full subcategory of
$\cd(C)$ whose objects are the simple perverse sheaves on $C$ which are equivariant for the
$G^0\T T\T T$ action 

(a) $(x,t,t'):(hU^*,h'U^*)\m(xht^nU^*,xh't'{}^nU^*)$
\nl
on $C$ (for some $n\in\NN^*_\kk$) or equivalently, are isomorphic to $\ucl^\sh_w[d_w]$ for
some $\cl\in\fs$ and some $w\in\WW$. Let $\cd^{cs}(C)$ be the subcategory of $\cd(C)$ whose
objects are those $K\in\cd(C)$ such that for any $j$, any simple subquotient of ${}^pH^jK$
is in $\hC$.

If $w,\cl$ are as above then $\ucl_w\in\cd^{cs}(C)$. Indeed this constructible sheaf is 
equivariant for the action (a) (for some $n$) hence so is each ${}^pH^j(\ucl_w)$.

We have a diagram $C\T C@<r<<(G^0/U^*)^3@>s>>C$ where 

$r(h_1U^*,h_2U^*,h_3U^*)=((h_1U^*,h_2U^*),(h_2U^*,h_3U^*))$,

$s(h_1U^*,h_2U^*,h_3U^*)=(h_1U^*,h_3U^*)$.
\nl
We define a bi-functor $\cd(C)\T\cd(C)@>>>\cd(C)$ by $A,A'\m A*A'=s_!r^*(A\bxt A')$. The 
"product" $A*A'$ is associative in an obvious sense. We show:

(b) $A,A'\m A*A'$ restricts to a bi-functor $\cd^{cs}(C)\T\cd^{cs}(C)@>>>\cd^{cs}(C)$.
\nl
Let $A,A'\in\cd^{cs}(C)$. To show that $A*A'\in\cd^{cs}(C)$ we may assume that 
$A,A'\in\hC$. Then each ${}^pH^j(A*A')$ is equivariant for the action (a) (for some $n$). 
This proves (b).

\subhead 40.5\endsubhead
For $w'\le w$ in $\WW$, $\l\in\ufs$, $\cl\in\l$ and $i\in\ZZ$ we show:

(a) {\it $\ch^i(\cl_w^\sh)|_{C_{w'}}\cong(\cl_{w'}(-i/2))^{\op N_{i,w',w,\l}}$}.
\nl
(Both sides are $0$ unless $i$ is even.) 

Let 

$\tC_w=\{(h,h')\in G^0\T G^0;h\i h'\in B^*\dw B^*\}\T\kk^*$, 

$\btC_w=\{(h,h')\in G^0\T G^0;h\i h'\in\ov{B^*\dw B^*}\}\T\kk^*$.
\nl
Now $\btC_w$ is an irreducible variety and we have a partition 
$\btC_w=\cup_{w';w'\le w}\tC_{w'}$ with $\tC_w$
smooth, open dense in $\btC_w$. Define $\bd:\btC_w@>>>\bC_w$, $d:\tC_w@>>>\bC_w$ by 
$(h,h',z)\m(hU^*,h'U^*)$. Let $\tcl_w=d^*\cl_w$, a local system on $\tC_w$. Let 
$\tcl_w^\sh=IC(\btC_w,\tcl_w)$. Since $d,\bd$ are principal $U^*\T\kk^*$-bundles it is 
enough to show

(b) {\it $\ch^i(\tcl_w^\sh)|_{\tC_{w'}}\cong(\tcl_{w'}(-i/2))^{\op N_{i,w',w,\l}}$}.
\nl
(Both sides are $0$ unless $i$ is even.) 

We choose $\k\in\Hom(T,\kk^*),\ce\in\fs(\kk^*)$ such that $\cl\cong\k^*\ce$, see 28.1(c).

Now $B^*$ acts on $(B^*\dw B^*)\T\kk^*$ and on $(\ov{B^*\dw B^*})\T\kk^*$ by 
$t_1u:(g,z)\m(g(t_1u)\i,\k(t_1)z)$ where $t_1\in T$, $u\in U^*$. Let 
$\bPP_w^\k=((\ov{B^*\dw B^*})\T\kk^*)/B^*$, $PP_w^\k=((B^*\dw B^*)\T\kk^*)/B^*$. Now 
$\PP_w^\k$ is a smooth open dense subvariety of the irreducible variety $\bPP_w^k$ and 
$\bPP_w^\k=\cup_{w';w'\le w}\PP_{w'}^\k$ is a partition. The morphism 
$(B^*\dw B^*)\T\kk^*@>>>\kk^*$ given by $(g,z)\m\k(\g_{\dw}(g))z$ factors through a 
morphism $\ph:\PP_w^\k@>>>\kk^*$. Let $\ce_w^\k=\ph^*\ce$, a local system of rank $1$ on 
$\PP_w^\k$. Let $\ce_w^{\k\sh}=IC(\bPP_w^\k,\ce_w^\k)\in\cd(\bPP_w^\k)$. From 
\cite{\CRG, 1.24} we see that

(c) {\it $\ch^i(\ce_w^{\k\sh})|_{\PP_{w'}^\k}\cong
(\ce_{w'}^\k(-i/2))^{\op N_{i,w',w,\l}}$}.
\nl
(Both sides are $0$ unless $i$ is even.) 

We can find $n\in\NN^*_\kk$ such that $\ce\in\fs_n(\kk^*)$. Define $\bc:\btC_w@>>>\bPP_w$,
$c:\btC_w@>>>\bPP_w$ by $(h,h',z)\m B^*-\text{orbit of }(h\i h',z^n)$. Now $\bc,c$ are 
locally trivial fibrations with smooth fibres of pure dimension. Hence (b) follows from (c)
provided that we can show that $c^*\ce_{w'}^\k=\tcl_{w'}$ for $w'\le w$. We may assume that
$w=w'$. We have a commutative diagram
$$\CD
\PP_w^\k@<c<<\tC_w\T\kk^*@>d>>C_w\\
@V\ph VV    @V\ph'VV   @V\k\ps VV\\
\kk^*@<c'<<\kk^*\T\kk^*@>d'>>\kk^* \endCD$$
with $\ph,\ps,c,d$ as above, $\ph'(h,h',z)=(\k(\g_{\dw}(h\i h')),z)$, $c'(z',z)=z'z^n$, 
$d'(z',z)=z'$. Using this and the definitions we have $\tcl_w=\ph'{}^*d'{}^*\ce$, 
$c^*\ce_w=\ph'{}^*c'{}^*\ce$. It remains to show that $d'{}^*\ce=c'{}^*\ce$. This expresses
the fact that $\ce$ is equivariant for the $\kk^*$-action $z_1:z\m z_1^nz$ on $\kk^*$ which
follows from $\ce\in\fs_n(\kk^*)$. This proves (b) hence (a).

\subhead 40.6\endsubhead
Let $w,w'\in\WW$, $\cl,\cl'\in\fs$. We set $L=\ucl_w*\ucl'_{w'}\in\cd^{cs}(C)$. Let
$$X=
\{(h_1U^*,h_2U^*,h_3U^*)\in(G^0/U^*)^3;h_1\i h_2\in B^*\dw B^*,h_2\i h_3\in B^*\dw'B^*\},$$
$$\align&\bX=\{(h_1U^*,h_2B^*,h_3U^*)\in G^0/U^*\T G^0/B^*\T G^0/U^*;\\&
h_1\i h_2\in B^*\dw B^*,h_2\i h_3\in B^*\dw'B^*\}.\endalign$$
We have a commutative diagram with a cartesian square 
$$\CD
X@>f>>\bX@>\bsi>>C\\
@V\t VV  @V\bat VV @.\\
T\T T@>f'>>T@.{}\endCD$$
where $f$ is given by $(h_1U^*,h_2U^*,h_3U^*)\m(h_1U^*,h_2B^*,h_3U^*)$, 

$f'$ is $(t,t')\m\Ad(\dw')\i(t)t'$,

$\t$ is $(h_1U^*,h_2U^*,h_3U^*)\m(t,t')$ with 
$h_1\i h_2\in U^*\dw tU^*,h_2\i h_3\in U^*\dw't'U^*$,

$\bat$ is $(h_1U^*,h_2B^*,h_3U^*)\m\Ad(\dw')\i(t)t'$ with $t,t'$ as in the definition of 
$\t$, 

$\bsi$ is $(h_1U^*,h_2B^*,h_3U^*)\m(h_1U^*,h_3U^*)$. 
\nl
From the definitions we have $L=\bsi_!f_!\t^*(\cl\bxt\cl')$. Using the diagram above, we 
have $L=\bsi_!\bat^*f'_!(\cl\bxt\cl')$. From the definitions we see that either (i) or (ii)
below holds:

(i) $\cl\not\cong(\Ad(\dw')\i)^*\cl'$ and $f'_!(\cl\bxt\cl')=0$;

(ii) $\cl\cong(\Ad(\dw')\i)^*\cl'$ and $\cl\bxt\cl'=f'{}^*\cl'$.
\nl
If (i) holds then $K=0$. If (ii) holds then, as in 32.16, we have
$$f'_!(\cl\bxt\cl')=f'_!f'{}^*\cl'
=\cl'\ot f'_!\bbq\Bpq\{\cl'\ot\ch^e(f'_!\bbq)[-e],e\in\ZZ\},$$
$$\cl'\ot\ch^e(f'_!\bbq)[-e]\Bpq\{\cl'(\rr-e),\do,\cl'(\rr-e), 
(\bin{\rr}{2\rr-e}\text{ copies})\}.$$
Setting $\bL=\bsi_!\bat^*(\cl')$, it follows that
$$L\Bpq\{\bL(\rr-e)[-e],\do,\bL(\rr-e)[-e],(\bin{\rr}{2\rr-e}\text{ copies}),e\in\ZZ\}.$$
We now consider $\bL$ for certain choices of $w,w'$.

If $w,w'$ satisfy $l(ww')=l(w)+l(w')$ then $\bsi$ restricts to an isomorphism
$\bX@>>>C_{ww'}$ and $\bL=\ucl'_{ww'}$. 

Now assume that $\a,\cha,s_\a$ are as in 28.3 and that $w=w'=s_\a\in\II$. We have

$\bL\Bpq\{j_{u!}\bL_u;u\in\WW\}$
\nl
where $j_u:C_u@>>>C$ is the inclusion and $\bL_u=j_u^*\bL$. Let $\bX_u=\bsi\i(C_u)$. Then 
$\bL_u=\bsi_{u!}\bat_u^*(\cl')$ where $\bsi_u:\bX_u@>>>C_u$, $\bat_u:\bX_u@>>>T$ are the 
restrictions of $\bsi,\bat$.

If $u\n\{1,s_\a\}$ then $\bX_u=\em$ and $\bL_u=0$. If $u=1$ then $\bsi_u:\bX_u@>>>C_u$ is 
an affine line bundle and $\bat_u^*(\cl')=\bsi_u^*\cl'_u$; hence 
$\bsi_{u!}\bat_u^*(\cl')=\bsi_{u!}\bsi_u^*\cl'_u=\cl'_u[[-1]]$. If $u=s_\a$ then 
$\bsi_u:\bX_u@>>>C_u$ is a principal $\kk^*$-bundle and either (iii) or (iv) below holds:

(iii) $\cha^*\cl'\not\cong\bbq$ and $\bsi_{u!}\bat_u^*(\cl')=0$,

(iv) $\cha^*\cl'\cong\bbq$ and $\bat_u^*(\cl')=\bsi_u^*\cl'_u$. 
\nl
If (iv) holds then, as in case (ii) above, we have

$\bsi_{u!}\bat_u^*(\cl')=\bsi_{u!}\bsi_u^*\cl'_u=\cl'_u\ot\bsi_{u!}\bbq\Bpq
\{\cl'_u\ot\ch^e(\bsi_{u!}\bbq)[-e],e\in\ZZ\}$,

$\cl'_u\ot\ch^e(\bsi_{u!}\bbq)[-e]\Bpq\{\cl'_u(1-e),\do,\cl'_u(1-e),(\bin{1}{2-e}
\text{ copies})\}$.

\subhead 40.7\endsubhead
In this subsection we assume that $\kk$ is an algebraic closure of a finite field. Now the
$\ca$-module $\fK(C)$ is defined as in 36.8 (the character sheaves on $C$ are taken to be 
the objects in $\hC$).

For $(w,\l)\in\WW\T\ufs$, let $[w;\l]$ be the basis element of $\fK(C)$ given by
$\ucl_w^\sh[[d_w/2]]$; we choose $\cl\in\l$ and we regard $\ucl_w,\ucl_w^\sh$ as mixed 
complexes on $C$ whose restriction to $C_w$ is pure of weight $0$; then 
$gr(\ucl_w),gr(\ucl_w^\sh)$ are defined in $\fK(C)$ as in 36.8.  We denote these elements 
of $\fK(C)$ by $[w;\l]',[w;\l]'{}^\sh$ respectively. From 40.5(a) we see that 

(a) $(-v)^{d_w}[w;\l]=[w;\l]'{}^\sh
=\sum_{w'\in\WW}\sum_{i\in 2\ZZ}N_{i,w',w,\l}v^i[w';\l]'$ in $\fK(C)$. 
\nl
where $N_{i,w',w,\l}$ is as in 40.2(e). 

Let $r,s$ be as in 40.4. By 40.4(b), $s_!r^*:\cd(C\T C)@>>>\cd(C)$ restricts to a functor 
$\cd^{cs}(C\T C)@>>>\cd^{cs}(C)$ where the character sheaves on $C\T C$ are by definition 
complexes of the form $A\bxt A'$ with $A\in\hC,A'\in\hC$. Hence the $\ca$-linear map 
$gr(s_!r^*):\fK(C\T C)@>>>\fK(C)$ or equivalently $\fK(C)\ot_\ca\fK(C)@>>>\fK(C)$ is well 
defined. (We have canonically $\fK(C\T C)=\fK(C)\ot_\ca\fK(C)$.) We write $\x*\x'$ instead
of $gr(s_!r^*)(\x\bxt\x')$ where $\x,\x'\in\fK(C)$. Note that $\x,\x'\m\x*\x'$ defines an
associative $\ca$-algebra structure on $\fK(C)$.

Let $w,w'\in\WW$, $\l,\l'\in\ufs$. From 40.6 we see that:

{\it if $w'\l'\ne\l$ then $[w;\l]'*[w';\l']'=0$ in $\fK(C)$;}

{\it if $w'\l'=\l$ and $l(ww')=l(w)+l(w')$ then $[w;\l]'*[w',\l']'=(v^2-1)^\rr[ww';\l']'$ 
in $\fK(C)$;}

{\it if $s\in\II$ and $s\l'=\l$ then 
$[s;\l]'*[s,\l']'=(v^2-1)^\rr(v^2[1;\l']'+(v^2-1)c[s;\l']')$ 
where $c=1$ for $s\in\WW_{\l'}$ and $c=0$ for $s\n\WW_{\l'}$.}
\nl
Using this and (a), 40.1(a), 40.2(e), we see that 

(b) {\it the unique $\ca$-linear isomorphism $\o:\fK(C)@>>>H$ ($H$ as in 40.1) given by 
$[w,\l]'\m v^{l(w)}\tT_w1_\l$ for $w\in\WW$, $\l\in\ufs$, satisfies
$\o([w,\l])=(-v)^{-d_w}v^{l(w)}c_{w,\l}$ for $w\in\WW$, $\l\in\ufs$ and
$\o(x*x')=(v^2-1)^\rr\o(x)\o(x')$ for any $x,x'\in\fK(C)$.}

\subhead 40.8\endsubhead
For $w,w'\in\WW$ and $\l,\l'\in\ufs$ we have 

$c_{w,\l}c_{w',\l'}=\sum_{y\in\WW,\nu\in\ufs}\g^{w,\l;w',\l'}_{y,\nu}c_{y,\l}$
\nl
in the algebra $H$. Here $\g^{w,\l;w',\l'}_{y,\nu}\in\ca$. We have:

(a) $\g^{w,\l;w',\l'}_{y,\nu}\in\NN[v,v\i]$.
\nl
By the arguments in 34.4-34.10 (with $D=G^0$) this is reduced to the analogous (well known)
statement for the structure constants of the algebra $H_\l^D$ with its basis $(c_w^\l)$ 
(see 34.2).

\subhead 40.9\endsubhead
For any $J\sub\II$ let $H_J$ be the $\ca$-submodule of $H$ spanned by 
$\{c_{w,\l};w\in\WW_J,\l\in\ufs\}$ or equivalently by $\{\tT_w1_\l;w\in\WW_J,\l\in\ufs\}$.
From the definitions we see that $H_J$ is a subalgebra of $H$. For any $J\sub\II,J'\sub\II$
we define a relation $\prq_{J,J'}$ on $\WW\T\ufs$ as follows. We say that 
$(y,\nu)\prq_{J,J'}(w,\l)$ if there exist $w_1\in\WW_J,w_2\in\WW_{J'}$, $\l_1,\l_2\in\ufs$
such that in the expansion (in the algebra $H$):

$c_{w_1,\l_1}c_{w,\l}c_{w_2,\l_2}=\sum_{y'\in\WW,\nu'\in\ufs}a_{y',\nu'}c_{y',\nu'}$
\nl
(with $a_{y',\nu'}\in\ca$) we have $a_{y,\nu}\ne0$.

Using the associativity of the product in $H$, the fact that $H_J,H_{J'}$ are subalgebras 
of $H$ and 40.8(a), we see that $\prq_{J,J'}$ is transitive. Using the formula 
$c_{1,w\l}c_{w,\l}c_{1,\l}=c_{w,\l}$ we see that it is reflexive. Thus, it is a preorder. 
Let $\si_{J,J'}$ be the equivalence relation attached to $\prq_{J,J'}$; thus, 
$(y,\nu)\si_{J,J'}(w,\l)$ if $(y,\nu)\prq_{J,J'}(w,\l)$ and $(w,\l)\prq_{J,J'}(y,\nu)$. The
equivalence classes for $\si_{J,J'}$ are called $(J,J')$-two-sided cells. The 
$(\II,\II)$-two sided cells in $\WW\T\ufs$ are also called two-sided cells.

\subhead 40.10\endsubhead
Let $w,w',w''\in\WW$, $\cl,\cl',\cl''\in\fs$. We set 
$K=\ucl_w*\ucl'_{w'}{}^\sh*\ucl''_{w''}\in\cd^{cs}(C)$. Let
$$\align&X=\{(h_1U^*,h_2U^*,h_3U^*,h_4U^*)\in(G^0/U^*)^4;\\&
h_1\i h_2\in B^*\dw B^*,h_2\i h_3\in\ov{B^*\dw'B^*},h_3\i h_4\in B^*\dw''B^*\},\endalign$$
an irreducible variety. Let $X_0$ be the smooth open dense subset of $X$ defined by the 
condition $h_2\i h_3\in B^*\dw'B^*$. Define $\s:X@>>>C$ by 

$(h_1U^*,h_2U^*,h_3U^*,h_4U^*)\m(h_1U^*,h_4U^*)$.
\nl
Define $\t:X_0@>>>T\T T\T T$ by 

$(h_1U^*,h_2U^*,h_3U^*,h_4U^*)\m(t,t',t'')$

with $h_1\i h_2\in U^*\dw tU^*,h_2\i h_3\in U^*\dw't'U^*,h_3\i h_4\in U^*\dw''t''U^*$.
\nl
Let $\cf=\t^*(\cl\bxt\cl'\bxt\cl'')$, a local system on $X_0$. Then 
$\cf^\sh:=IC(X,\cf)\in\cd(X)$ is defined and we have $K=\s_!\cf^\sh$.

Let $\bX$ (resp. $\bX_0$) be the the variety of all 
$(h_1U^*,h_2B^*,h_3B^*,h_4U^*)\in G^0/U^*\T G^0/B^*\T G^0/B^*\T G^0/U^*$ that satisfy the 
same equations as those defining $X$ (resp. $X_0$). Note that $\bX$ is irreducible and
$\bX_0$ is an open dense smooth subset of $\bX$. We have a cartesian diagram
$$\CD
X@>f>>\bX@>\bsi>>C\\
@AAA   @AAA  @.\\
X_0@>f_0>>\bX_0@.{}\\
@V\t VV  @V\bat VV @.\\
T\T T\T T@>f'>>T@.{}\endCD$$
where $X_0@>>>X,\bX_0@>>>\bX$ are the obvious imbeddings,

$f,f_0$ are given by $(h_1U^*,h_2U^*,h_3U^*,h_4U^*)\m(h_1U^*,h_2B^*,h_3B^*,h_4U^*)$, 

$f'$ is $(t,t',t'')\m\Ad(\dw'\dw'')\i(t)\Ad(\dw'')\i(t')t''$,

$\bat$ is $(h_1U^*,h_2B^*,h_3B^*,h_4U^*)\m\Ad(\dw'\dw'')\i(t)\Ad(\dw'')\i(t')t''$ with 
$t,t',t''$ as in the definition of $\t$, 

$\bsi$ is $(h_1U^*,h_2B^*,h_3B^*,h_4U^*)\m(h_1U^*,h_4U^*)$. 
\nl
Assume that $\cl\cong(\Ad(\dw')\i)^*\cl'$ and $\cl'\cong(\Ad(\dw'')\i)^*\cl''$. Then 
$\cl\bxt\cl'\bxt\cl''=f'{}^*\cl''$. We have $\cf=\t^*f'{}^*\cl''=f_0^*\bat^*\cl''$. Since 
$f$ is a principal $T\T T$-bundle and $X_0=f\i(\bX_0)$ it follows that
$\cf^\sh=f^*IC(\bX,\bat^*\cl'')$. Note that 
$f_!\bbq\Bpq\{\ch^e(f_!\bbq)[-e],2\rr\le e\le4\rr\},$
$$\ch^e(f_!\bbq)\Bpq\{\bbq(2\rr-e),\do,\bbq(2\rr-e),(\bin{2\rr}{4\rr-e}\text{ copies})\}.$$
Hence setting $\bK=\bsi_!(IC(\bX,\bat^*\cl''))$ we have 
$$K=\s_!f^*IC(\bX,\bat^*\cl'')=\bsi_!f_!f^*IC(\bX,\bat^*\cl'')
=\bsi_!(IC(\bX,\bat^*\cl'')\ot f_!\bbq),$$
$$K\Bpq\{\bK(2\rr-e)[-e],\do,\bK(2\rr-e)[-e],(\bin{2\rr}{4\rr-e}\text{ copies}),
2\rr\le e\le4\rr\}.\tag a$$
We now show:

(b) {\it if $A\in\hC$ is such that $A\dsv\bK$, then $A\dsv K$.}
\nl
We may regard $\cl,\cl',\cl''$ as mixed local systems (with respect to a rational structure
over a sufficiently large finite subfield of $\kk$) which are pure of weight $0$. Then 
$K,\bK$ are naturally mixed complexes and (a) is compatible with the mixed structures. For
any mixed perverse sheaf $P$, let $P_h$ be the subquotient of $P$ of pure weight $h$. We 
can find $h\in\ZZ$ such that $A\dsv{}^pH^j(\bK)_h$ for some $j\in\ZZ$; moreover we may 
assume that $h$ is maximum possible. Note that 
$A\dsv{}^pH^{j+4\rr}(\bK[-4\rr](-2\rr))_{h+2\rr}$ and
$A\not\dsv{}^pH^{j'}(\bK[-e](2\rr-e))_{h+2\rr}$ for $2\rr\le e<4\rr$ and any $j'$; hence
from (a) we see that $A\dsv{}^pH^{j+4\rr}(K)_{h+2\rr}$. In particular, $A\dsv K$, and (b) 
is proved.

\subhead 40.11\endsubhead
Let $w,w'\cl,\cl',X,\bX,\t$ be as in 40.6. We set 
$\LL=\ucl_w^\sh*\ucl'_{w'}{}^\sh\in\cd^{cs}(C)$. Let $A=\ucl''_{w''}{}^\sh[d_{w''}]$. We 
show:

(a) {\it If $A\dsv\LL$ then $[w'',\l'']$ appears with non-zero coefficient in the 
expansion of the product $[w,\l]*[w',\l']$ in terms of the basis $([y,\nu])$ of $\fK(C)$.}
\nl
Let
$$\XX=\{(h_1U^*,h_2U^*,h_3U^*)\in(G^0/U^*)^3;h_1\i h_2\in\ov{B^*\dw B^*},
h_2\i h_3\in\ov{B^*\dw'B^*}\},$$
$$\align&\bar{\XX}=\{(h_1U^*,h_2B^*,h_3U^*)\in G^0/U^*\T G^0/B^*\T G^0/U^*;\\&
h_1\i h_2\in\ov{B^*\dw B^*},h_2\i h_3\in\ov{B^*\dw'B^*}\}.\endalign$$
Note that $X$ (resp. $\bX$) is naturally an open dense subset of $\XX$ (resp. $\bar{\XX}$).
Define $\s':\XX@>>>C$ by $(h_1U^*,h_2U^*,h_3U^*)\m(h_1U^*,h_3U^*)$. Define 
$\bsi':\bar{\XX}@>>>C$ by $(h_1U^*,h_2B^*,h_3U^*)\m(h_1U^*,h_3U^*)$. Let 
$\cf=\t^*(\cl\bxt\cl')$, a local system on $X$. Then $\cf^\sh:=IC(\XX,\cf)\in\cd(\XX)$ is 
defined and we have $\LL=\s'_!\cf^\sh$. We have a cartesian diagram
$$\CD
\XX@>\tf>>\bar{\XX}@>\bsi'>>C\\
@AAA   @AAA  @.\\
X@>f>>\bX@.{}\\
@V\t VV  @V\bat VV @.\\
T\T T@>f'>>T@.{}\endCD$$
where $X@>>>\XX,\bX@>>>\bar{\XX}$ are the obvious imbeddings,
$f,f',\bat$ are as in 40.6 and $\tf$ is the obvious map.

Assume first that 40.6(i) holds. Let $m':T\T\XX@>>>\XX$ be the free $T$-action 
$t_1:(h_1U^*,h_2U^*,h_3U^*)\m(h_1U^*,h_2t_1\i U^*,h_3U^*)$. This restricts to a free 
$T$-action $m:T\T X@>>>X$. Define a free $T$ action $m_0:T\T(T\T T)@>>>T\T T$ by
$t_1:(t,t')\m(t_1\i t,\Ad(\dw')\i(t_1)t'$. Then $m,m_0$ are compatible with $\t$. By our 
assumption we have $m_0^*(\cl\bxt\cl')=\cl_0\bxt\cl\bxt\cl'$ where $\cl_0\in\fs(T)$,
$\cl_0\not\cong\bbq$. It follows that $m^*(\cf)\cong\cl_0\bxt\cf$. From the properties of
intersection cohomology we then have $m'{}^*(\cf^\sh)\cong\cl_0\bxt\cf^\sh$. Let 
$r:T\T\XX@>>>\XX$ be the second projection. Since $\cl_0\in\fs(T)$, $\cl_0\not\cong\bbq$,
we have $r_!(\cl_0\bxt\cf^\sh)=0$. Hence $r_!m'{}^*(\cf^\sh)=0$. Since $m',f',r,f'$ form a
cartesian diagram we must have $f'{}^*f'(\cf^\sh)=0$. Since $f'$ is a principal $T$-bundle
we deduce that $f'_!(\cf^\sh)=0$. We have $\LL=\bsi'_!f'_!(\cf^\sh)$ hence $\LL=0$. In this
case (a) is clear.

Assume next that 40.6(ii) holds. Then $\cl\bxt\cl'=f'{}^*\cl'$ and 
$\cf=\t^*f'{}^*\cl'=f^*\bat^*\cl'$. Since $f'$ is a principal $T$-bundle and 
$X=f'{}\i(\bX)$ it follows that $\cf^\sh=f'{}^*IC(\bar{\XX},\bat^*\cl')$. Note that 
$f'_!\bbq\Bpq\{\ch^e(f'_!\bbq)[-e],\rr\le e\le2\rr\},$
$$\ch^e(f'_!\bbq)\Bpq\{\bbq(\rr-e),\do,\bbq(\rr-e),(\bin{\rr}{2\rr-e}\text{ copies})\}.$$
Hence setting $\bar{\LL}=\bsi'_!(IC(\bar{\XX},\bat^*\cl'))$ we have 
$$\LL=\s'_!f'{}^*IC(\bar{\XX},\bat^*\cl')=\bsi'_!f'_!f'{}^*IC(\bar{\XX},\bat^*\cl')
=\bsi'_!(IC(\bar{\XX},\bat^*\cl')\ot f'_!\bbq),$$
$$\LL\Bpq\{\bar{\LL}(\rr-e)[-e],\do,\bar{\LL}(\rr-e)[-e],(\bin{\rr}{2\rr-e}\text{ copies}),
\rr\le e\le2\rr\}.$$
Since $A\dsv\LL$, this shows that $A\dsv\bar{\LL}$. We regard $\cl'$ as a pure local
system of weight $0$. Then $\bar{\LL}=\bsi'_!(IC(\bar{\XX},\bat^*\cl'))$ is again pure of
weight $0$, since $\bsi'$ is proper (see \cite{\BBD}). Hence the coefficient with which $A$
appears in the expansion of $gr(\bar{\LL})$ is a polynomial in $-v$ with coefficients given
by the multiplicities of $A$ in the various ${}^pH^j(\bar{\LL})$; in particular, $A$ 
appears with coefficient $\ne 0$ in $gr(\bar{\LL})$. On the other hand the arguments above
show that $[w,\l]*[w',\l']=(v^2-1)^\rr gr(\bar{\LL})$. It follows that $A$ appears with 
coefficient $\ne 0$ in $[w,\l]*[w',\l']$. This proves (a).

\head 41. Character sheaves and two-sided cells\endhead 
\subhead 41.1\endsubhead
In this section we preserve the notation of 40.3. We fix a connected component $D$ of $G$ 
and we pick $\d\in N_DB^*\cap N_DT$. We write $\e$ instead of $\e_D:\WW@>>>\WW$. For 
$w\in\WW$ we set

$Z^w_{\em,D}=\{(B,B',xU_B)\in Z_{\em,D};\po(B,B')=w\}$.
\nl
(This is the same as ${}^{w\i}Z_{\em,D}$ in 36.2.) Define $\x_D:C@>>>Z_{\em,D}$ by 
$(hU^*,h'U^*)\m(hB^*h\i,h'B^*h'{}\i,h'\d h\i U_{hB^*h\i})$, a principal $T$-bundle for the
free $T$-action on $C$ given by $t:(hU^*,h'U^*)@>>>(htU^*,h'(\d t\d\i)U^*)$. 

Since $\x_D\i(Z^w_{\em,D})=C_w$, $\x_D$ restricts to a principal $T$-bundle 
$\x_{D,w}:C_w@>>>Z^w_{\em,D}$. We have a commutative diagram
$$\CD
T@<\ps<<C_w@>=>>                     C_w\\
@V\z VV   @Aj'AA                                 @V\x_{D,w}VV\\
\dd@<pr_2<<G^0/(U^*\cap\dw U^*\dw\i)\T\dd@>j>>Z^w_{\em,D}\endCD$$
where $\ps$ is as in 40.3,

$\dd=\dw\d T$,

$j(f(U^*\cap\dw U^*\dw\i),s)=(fB^*f\i,f\dw B^*\dw\i f\i,fsf\i U_{fB^*f\i})$,

$j'(f(U^*\cap\dw U^*\dw\i),s)=(fU^*,fs\d\i U^*)$,

$\z(t)=\dw\d(\d\i t\d)$.   
\nl
Note that the lower row in the diagram is as in 36.2(a).

Define $\io:\dd@>>>T$ by $\io(\dw\d t)=t$ where $t\in T$. If $\cl\in\fs$ is such that
$\Ad((\dw d)\i)^*\cl\cong\cl$ then $pr_2^*\io^*(\cl)$ is a local system on
$G^0/(U^*\cap\dw U^*\dw\i)\T\dd$, equivariant for the $T$-action
$t_0:(f(U^*\cap\dw U^*\dw\i),s)=(ft_0\i(U^*\cap\dw U^*\dw\i),t_0st_0\i)$ on
$G^0/(U^*\cap\dw U^*\dw\i)\T\dd$, which makes $j$ a principal $T$-bundle. It follows that
there is a well defined local system $\dcl_w$ (of rank $1$) on $Z^w_{\em,D}$ such that
$j^*\dcl_w=pr_2^*\io^*(\cl)$. We show:

(a) $\x_{D,w}^*(\dcl_w)=(\Ad(\d\i)^*\cl)_w$.
\nl
Since $j'$ is an isomorphism it is enough to show that
$j'{}^*\x_{D,w}^*(\dcl_w)=j'{}^*(\Ad(\d\i)^*\cl)_w$ or that
$j^*\dcl_w=j'{}^*(\Ad(\d\i)^*\cl)_w$ or that $pr_2^*\io^*\cl=j'{}^*\ps^*(\Ad(\d\i)^*)\cl)$
or that $j'{}^*\ps^*\z^*\io^*\cl=j'{}^*\ps^*(\Ad(\d\i)^*)\cl)$. It is enough to show that 
$\z^*\io^*\cl=\Ad(\d\i)^*\cl$. This follows from $\Ad(\d\i)=\io\z:T@>>>T$.

Let $h_w:Z^w_{\em,D}@>>>Z_{\em,D}$, $\bh_w:\bZ^w_{\em,D}@>>>Z_{\em,D}$ be the inclusions
($\bZ^w_{\em,D}=\cup_{w';w'\le w}Z^{w'}_{\em,D}$ is the closure of $Z^w_{\em,D}$ in 
$Z_{\em,D}$.) Let $\udcl_w=h_{w!}\dcl_w,\udcl^\sh_w=\bh_{w!}\dcl^\sh_w$. Using (a) and the
fact that $\x_D$ is a principal $T$-bundle we deduce 

(b) $\x_D^*(\udcl_w)=\un{(\Ad(\d\i)^*\cl)}_w$,

(c) $\x_D^*(\udcl^\sh_w)=\un{(\Ad(\d\i)^*\cl)}^\sh_w$.
\nl
Now let $D'$ be another connected component of $G$. We pick $\d'\in N_{D'}B^*\cap N_{D'}T$.
We have a commutative diagram with a cartesian right square
$$\CD
C\T C@<r<<(G^0/U^*)^3@>s>>C\\ 
@V\x_D\T\x_{D'}VV               @V\x_0VV              @V\x_{D'D}VV\\
Z_{\em,D}\T Z_{\em,D'}@<b_1<<Z_0@>b_2>>Z_{\em,D'D}   \endCD$$
where $r,s$ are as in 40.4, $Z_0,b_1,b_2$ are as in 32.5 (with $J=\em$) and 
$$\align&\x_0(h_1U^*,h_2U^*,h_3U^*)\\&
=(h_1B^*h_1\i,h_2B^*h_2\i,h_3B^*h_3\i,
h_2\d h_1\i U_{h_1B^*h_1\i},h_3\d'h_2\i U_{h_2B^*h_2\i}).\endalign$$
Hence, if $A\in\cd(Z_{\em,D}),A'\in\cd(Z_{\em,D'})$, then 
$\x_{D'D}^*b_{2!}b_1^*(A\bxt A')=s_!r^*(\x_D^*A\bxt\x_{D'}^*A')$, or equivalently

(d) $\x_{D'D}^*(A*A')=(\x_D^*A)*(\x_{D'}^*A')$.

\subhead 41.2\endsubhead
Let $u\in\WW$. Let
$$\align\Up_u=&\{(B,B',g(U_B\cap U_{B'});\\&
B\in\cb,B'\in\cb,g(U_B\cap U_{B'})\in D/(U_B\cap U_{B'}),\po(B,B')=u\}\endalign$$
and let $\Ph_u:\cd(Z_{\em,D})@>>>\cd(Z_{\em,D})$ be the composition $\fh_!\fj^*$ where
$\fj:\Up_u@>>>Z_{\em,D}$ is $(B,B',g(U_B\cap U_{B'})\m(B,gBg\i,gU_B)$ and
$\fh:\Up_u@>>>Z_{\em,D}$ is 

$(B,B',g(U_B\cap U_{B'})\m(B',gB'g\i,gU_{B'})$.
\nl
(A special case of definitions in 37.1.) Let
$$\align\Up'=&\{(B',B,\tB,\tB',gU_{B'});
B'\in\cb,B\in\cb,\tB\in\cb,\tB'\in\cb,\\&
gU_{B'}\in D/U_{B'},\po(B',B)=u\i,\po(\tB,\tB')=\e(u),gB'g\i=\tB'\},\endalign$$
$$s:\Up_u@>>>\Up',(B,B',g(U_B\cap U_{B'})\m(B',B,gBg\i,gB'g\i,gU_{B'}).$$
Note that $s$ is an isomorphism. (We show this only at the level of sets. Define 
$s':\Up'@>>>\Up_u$ by $(B',B,\tB,\tB',gU_{B'})\m(B,B',x(U_B\cap U_{B'}))$ where $x\in D$ is
such that $xBx\i=\tB$, $xU_{B'}=gU_{B'}$. This is well defined and clearly an inverse of
$s$.) It follows that $\fh_!\fj^*=\fh'_!\fj'{}^*$ where 

$\fh'=\fh s':\Up'@>>>Z_{\em,D}$ is $(B',B,\tB,\tB',gU_{B'})\m(B',\tB',gU_{B'})$, 

$\fj'=\fj s':\Up'@>>>Z_{\em,D}$ is $(B',B,\tB,\tB',gU_{B'})\m(B,\tB,xU_B)$
\nl
and $x\in D$ is such that 
$xBx\i=\tB,xU_{B'}=gU_{B'}$ (then $x(U_B\cap U_{B'})$ is well defined). We have a
commutative diagram with a cartesian right square
$$\CD
C   @<\tj<<    \tC@>\tih>>      C\\
@V\x_DVV      @V\x'VV    @V\x_DVV\\
Z_{\em,D}@<\fj'<< \Up'_u@>\fh'>>Z_{\em,D}\endCD$$
where $\x_D$ is as in 41.1,
$$\align\tC=&\{(h_1U^*,h_2B^*,h_3B^*,h_4U^*)\in(G^0/U^*)^4;\\&
 h_1\i h_2\in B^*\du\i B^*,h_3\i h_4\in B^*\d\du\d\i B^*\},\endalign$$
$\tih$ is $(h_1U^*,h_2B^*,h_3B^*,h_4U^*)\m(h_1U^*,h_4U^*)$,

$\x'$ is 
$$\align&(h_1U^*,h_2B^*,h_3B^*,h_4U^*)\\&
\m(h_1B^*h_1\i,h_2B^*h_2\i,h_3B^*h_3\i,h_4B^*h_4\i,h_4\d h_1\i U_{h_1B^*h_1\i}),\endalign$$
$\tj$ is $(h_1U^*,h_2B^*,h_3B^*,h_4U^*)\m(h_2t\i U^*,h_3\tit U^*)$ 

where $t,\tit\in T$ are
given by $h_1\i h_2\in U^*\du\i tU^*$, $h_3\i h_4\in U^*\tit\d\du\d\i U^*$.
\nl
We see that for $A\in\cd(Z_{\em,D})$ we have
$$\x_D^*\Ph_u(A)=\x_D\fh_!\fj^*A=\x_D^*\fh'_!\fj'{}^*A=\tih_!\x'{}^*\fj'{}^*A=
\tih_!\tj^*\x_D^*A.$$
Taking here $A=\udcl_w^\sh$ (with $w\in\WW,\l\in\ufs,\cl\in\l$ with $w\uD\l=\l$) and using
41.1(c) we obtain $\x_D^*\Ph_u(\udcl_w^\sh)=\tih_!\tj^*(\un{(\Ad(\d\i)^*\cl)}_w^\sh)$ or 
equivalently $\x_D^*\Ph_u(\udcl_w^\sh)=\bsi_!\tj'{}^*((\Ad(\d\i)^*\cl)_w^\sh)$ where

$\bX=\{(h_1U^*,h_2B^*,h_3B^*,h_4U^*)\in\tC;h_2\i h_3\in\ov{B^*\dw B^*}\}$
\nl
and $\tj':\bX@>>>\bC_w,\bsi:\bX@>>>C$ are the restrictions of $\tj,\tih$. Let 

$\bX_0=\{(h_1U^*,h_2B^*,h_3B^*,h_4U^*)\in\tC;h_2\i h_3\in B^*\dw B^*\}$
\nl
and let
$\tj'_0:\bX_0@>>>C_w$ be the restriction of $\tj$. Let $\cf_0=\tj'_0{}^*(\Ad(\d\i)^*\cl)$,
a local system on $\bX_0$. Since $\tj'$ is a fibration with smooth connected fibres we have
$\tj'{}^*((\Ad(\d\i)^*\cl)_w^\sh)=IC(\bX,\cf_0)$. Thus, 
$\x_D^*\Ph_u(\udcl_w^\sh)=\bsi_!(IC(\bX,\cf_0))$. From the definitions we see that 
$\cf_0=\bat^*\cl''$ hence $\bsi_!(IC(\bX,\cf_0))=\bK$ and 
$$\x_D^*\Ph_u(\udcl_w^\sh)=\bK\tag a$$ 
where $\bat^*\cl',\bK$ are given as in 40.10 in terms of 

$(u\i,\cl),(w,\Ad(\d\i)^*\cl),(\e(u),\Ad(\d\du\d\i)^*\Ad(\d\i)^*\cl)$
\nl
instead of $(w,\cl),(w',\cl'),(w'',\cl'')$. 

\subhead 41.3\endsubhead
For $J\sub\II$ let $\cd^{cs}_J(C)$ be the subcategory of $\cd^{cs}(C)$ whose objects are 
those $K\in\cd(C)$ such that for any $j$, any simple subquotient of ${}^pH^jK$ is 
isomorphic to $\ucl_w^\sh$ for some $\cl\in\fs$ and some $w\in\WW_J$. 

Let $J,J'\sub\II$. Let $K\in\cd^{cs}_J(C),K'\in\cd^{cs}_{J'}(C)$, and let $w',w''\in\WW$, 
$\l',\l''\in\ufs$, $\cl'\in\l',\cl''\in\l''$. Let $A=\ucl''_{w''}{}^\sh[d_{w''}]$. We show:

(a) {\it If (i) $A\dsv K*\ucl'_{w'}{}^\sh[d_{w'}]$ or (ii) 
$A\dsv\ucl'_{w'}{}^\sh[d_{w'}]*K'$ or (iii) $A\dsv K*\ucl'_{w'}{}^\sh[d_{w'}]*K'$ then 
$(w'',\l'')\prq_{J,J'}(w',\l')$.}
\nl
For the proof we may assume that $\kk$ is an algebraic closure of a finite field. Then the
results in 40.7 are applicable. We first consider the case (i). In this case we can find 
$\cl\in\fs,w\in\WW_J$ such that $A\dsv\ucl_w^\sh[d_w]*\ucl'_{w'}{}^\sh[d_{w'}]$. By 
40.11(a), $[w'',\l'']$ appears with non-zero coefficient in the expansion of the product 
$[w,\l]*[w',\l']$ in terms of the basis $([y,\nu])$ of $\fK(C)$. Applying $\o$ (see 
40.7(b)) we see that $c_{w'',\l''}$ appears with non-zero coefficient in the expansion of 
the product $c_{w,\l}c_{w',\l'}$ in terms of the basis $(c_{y,\nu})$ of $H$ and the desired
result follows. Case (ii) is treated in an entirely similar way. We now consider case 
(iii). In this case we must have $A\dsv A'*K'$ for some simple perverse sheaf $A'$ such 
that $A'\dsv K*\ucl'_{w'}{}^\sh[d_{w'}]$. We have $A'=\ucm_y^\sh[d_y]$ where $y\in\WW$, 
$\cm\in\fs$. Let $\nu$ be the isomorphism class of $\cm$. From case (ii) applied to 
$A\dsv A'*K'$ we see that $(w'',\l'')\prq_{J,J'}(y,\nu)$. From case (i) applied to 
$A'\dsv K*\ucl'_{w'}{}^\sh[d_{w'}]$ we see that $(y,\nu)\prq_{J,J'}(w',\l')$. Combining 
these two inequalities we obtain $(w'',\l'')\prq_{J,J'}(w',\l')$, as desired.

\subhead 41.4\endsubhead
Let $J\sub\II$. In the remainder of this section we write $\ff,\fe$ instead of
$\ff_{\em,J}:\cd(Z_{\em,D}@>>>\cd(Z_{J,D})$, $\fe_{\em,J}:\cd(Z_{J,D}@>>>\cd(Z_{\em,D})$. 
We note:

(a) {\it If $A\in\cd(Z_{J,D})$ then $\ff\fe(A)\cong A[m]\op A'$ for some $m\in\ZZ$ and some
$A'\in\cd(Z_{J,D})$.}
\nl
See \cite{\GI}, \cite{\MV} for the special case $D=G^0,J=\II$ and \cite{\PCS, 6.6} for the
general case. We show:

(b) {\it Let $A$ be a simple perverse sheaf on $Z_{J,D}$. Then $A\dsv\ff({}^pH^j(\fe(A)))$ 
for some $j\in\ZZ$.}
\nl
Assume that this is not true. As in \cite{\BBD, p.142}, for any $n\in\ZZ$ we have a 
distinguished triangle $({}^p\t_{\le n-1}\fe A,{}^p\t_{\le n}\fe A,{}^pH^n(\fe A)[-n])$ 
hence a distinguished triangle

$(\ff({}^p\t_{\le n-1}\fe A),\ff({}^p\t_{\le n}\fe A),\ff({}^pH^n(\fe A))[-n])$.
\nl
Using our assumption, we see that $A\dsv\ff({}^p\t_{\le n-1}\fe A)$ if and only if
$A\dsv\ff({}^p\t_{\le n}\fe A)$. Thus we have $A\dsv\ff({}^p\t_{\le n}\fe A)$ for some $n$
if and only if $A\dsv\ff({}^p\t_{\le n}\fe A)$ for any $n$. Since ${}^p\t_{\le n}\fe A=0$ 
for some $n$ we see that $A\ndsv\ff({}^p\t_{\le n}\fe A)$ for any $n$. Since
${}^p\t_{\le n}\fe A=\fe A$ for some $n$ we deduce that $A\ndsv\ff\fe A$. This contradicts
(a); (b) is proved.

We show:

(c) {\it If $A$ is a simple perverse sheaf on $Z_{J,D}$ then there exists a simple perverse
sheaf $A'$ on $Z_{\em,D}$ such that $A\dsv\ff(A'),A'\dsv\fe(A)$.}
\nl
By (b) we can find $i,j\in\ZZ$ such that $A\dsv{}^pH^i(\ff(P))$ where $P={}^pH^j(\fe(A))$.

Assume that $A\ndsv{}^pH^i(\ff(A'))$ for any simple subquotient $A'$ of $P$. We claim that
$A\ndsv{}^pH^i(\ff(P'))$ for any subobject $P'$ of $P$. We argue by induction on the length
of $P'$. If $P'$ has length $1$ the claim holds by assumption. If $P'$ has length $\ge 2$,
we can find a simple subobject $P''$ of $P'$. We have a distinguished triangle 
$(\ff(P''),\ff(P'),\ff(P'/P''))$. Hence we have an exact sequence 
${}^pH^i(\ff(P''))@>>>{}^pH^i(\ff(P'))@>>>{}^pH^i(\ff(P'/P''))$. By the induction
hypothesis, we have $A\ndsv{}^pH^i(\ff(P'')), A\ndsv{}^pH^i(\ff(P'/P''))$. Hence 
$A\ndsv{}^pH^i(\ff(P'))$. This proves the claim. In particular, $A\ndsv{}^pH^i(\ff(P))$,
contradicting the definition of $i,P$. 

We see that there exists a simple subquotient $A'$ of $P$ such that 
$A\dsv{}^pH^i(\ff(A'))$. Then $A'$ is as required by (c).

Let $\bd_w=\dim Z_{\em,D}^w$. Let 

(d) $A'=\udcl_w^\sh[\bd_w],A''=\udcm_y^\sh[\bd_y]\in\hZ_{\em,D}$, $\cl\in\l,\cm\in\nu$.
\nl
Here $w\uD\l=\l,y\uD\nu=\nu$. We show:

(e) {\it Let $A$ be a character sheaf on $Z_{J,D}$ such that $A\dsv\ff(A'),A''\dsv\fe(A)$.
Then $(y,\uD\nu)\prq_{J,J'}(w,\uD\l)$.}  
\nl
Since $\ff$ is proper, $\ff(A')$ is a semisimple complex (see \cite{\BBD}). Hence 
$\ff(A')\cong A[m]\op A_1$ for some $m\in\ZZ,A'\in\cd(Z_{J,D})$ and 
$\fe\ff(A')\cong\fe(A)[m]\op\fe(A_1)$. Hence from $A''\dsv\fe(A)$ we can deduce 
$A''\dsv\fe\ff(A')$. By 37.2 we have $\fe\ff(A')\Bpq\{\Ph_u(A')[[-m_u]];u\in\WW_J\}$ where
$m_u$ are certain integers. Hence for some $u\in\WW_J$ we have $A''\dsv\Ph_u(A')[[-m_u]]$ 
that is, $A''\dsv\Ph_u(A')$ and $\x_D^*A''[\rr]\dsv\x_D^*\Ph_u(A')[\rr]$. Hence using 
41.2(a) we have $\x_D^*A''[\rr]\dsv\bK$ where $\bK$ is as in the end of 41.2. Thus, 
$\ucm_y^\sh[d_y]\dsv\bK$. Using 40.10(b) we deduce that 
$$\ucm_y^\sh[d_y]\dsv\un{\Ad(\dw)\i)^*\Ad(\d\i)^*\cl}_{u\i}*
\un{(\Ad(\d\i)^*\cl)}_w^\sh*\un{\Ad(\d\du\d\i)^*\Ad(\d\i)^*\cl}_{\e(u)}.$$
Using this and 41.3(a) we see that (e) holds. 

We show:

(f) {\it Let $A$ be a character sheaf on $Z_{J,D}$. In the setup of (d) assume that
$A\dsv\ff(A')$, $A'\dsv\fe(A)$, $A\dsv\ff(A'')$, $A''\dsv\fe(A)$. Then
$(y,\uD\nu)\si_{J,J'}(w,\uD\l)$.}  
\nl
Applying (e) to $A',A''$ we see that $(y,\uD\nu)\prq_{J,J'}(w,\uD\l)$. Applying (e) to
$A'',A'$ (instead of $A',A''$) we see that $(w,\uD\l)\prq_{J,J'}(y,\uD\nu)$. Hence (f)
holds.

\mpb

From (c),(f) we see that there is a well defined map $A\m\boc_A$ from the set of character 
sheaves on $Z_{J,D}$ (up to isomorphism) to the set of $(J,J')$-two-sided cells in
$\WW\T\ucf$ where $\boc_A$ is the unique $(J,J')$-two-sided cell that contains 
$$\{(w,\uD\l)\in\WW\T\ufs;w\uD\l=\l,A\dsv\ff(\udcl_w^\sh[\bd_w]),\udcl_w^\sh[\bd_w]\dsv A\}
$$
(a non-empty set); here $\cl\in\l$.

\subhead 41.5\endsubhead
In the setup of 41.4, let $A$ be a character sheaf on $Z_{J,D}$. We show:

(a) {\it There exists $(w,\uD\l)\in\boc_A$ such that $w\uD\l=\l$,
$A\dsv\ff(\udcl_w^\sh[\bd_w])$. If $(w',\uD\l')\in\WW\T\ufs$ is such that $w'\uD\l'=\l'$, 
$A\dsv\ff(\udcl'_{w'}{}^\sh[\bd_{w'}])$ then $(w,\uD\l)\prq_{J,J'}(w',\uD\l')$. Here 
$\cl\in\l,\cl'\in\l'$.}

(b) {\it There exists $(w,\uD\l)\in\boc_A$ such that $w\uD\l=\l$,
$\udcl_w^\sh[\bd_w]\dsv\fe(A)$. If $(w',\uD\l')\in\WW\T\ufs$ is such that $w'\uD\l'=\l'$,
$\udcl'_{w'}{}^\sh[\bd_{w'}]\dsv\fe(A)$ then $(w',\uD\l')\prq_{J,J'}(w,\uD\l)$. Here 
$\cl\in\l,\cl'\in\l'$.}
\nl
Note that (a) follows immediately from 41.4(c),(e) and the definition of $\boc_A$. 
Similarly, (b) follows from 41.4(c),(e) and the definition of $\boc_A$. 

\subhead 41.6\endsubhead
In this subsection we assume that $J=\II$. The $\ca$ linear map $H@>>>H$ given by 

(a) $\tT_w1_\l\m\tT_{\e(w)}1_{\uD\l}$ for $w\in\WW,\l\in\ufs$
\nl
is an $\ca$-algebra isomorphism. It carries $c_{w,\l}$ to $c_{\e(w),\uD\l}$ for any 
$w\in\WW,\l\in\ufs$. It induces a bijection $\boc\m\boc'$ from the set of two-sided cells
in $\WW\T\ufs$ onto itself. We show:

(b) {\it If $A$ is a character sheaf on $D$ then $(\boc_A)'=\boc_A$.}
\nl
Consider the automorphism $\Ad(\d):D@>>>D$. From the definitions we see that for 
$(w,\l)\in\WW\T\ufs$ such that $w\uD\l=\l$ we have $A\dsv\ff(\udcl_w^\sh[\bd_w])$ if and 
only if $\Ad(\d\i)^*A\dsv\ff(\un{\dot{\Ad(\uD\i)^*\cl}}_{\e(w)}^\sh[\bd_w])$. Using this
and 41.5(a) we see that 

$\boc_{\Ad(\d\i)^*A}=(\boc_A)'$. 
\nl
It is then enough to show that 
$\Ad(\d\i)^*A\cong A$. By the $G^0$-equivariance of $A$ we have $m^*A\cong q^*A$ where 
$m:G^0\T D@>>>D$ is $(x,g)\m xgx\i$ and $q:G^0\T D@>>>D$ is $(x,g)\m g$. Define 
$r:D@>>>G^0\T D$ by $r(g)=(\d g\i,g)$. Then $r^*m^*A\cong r^*q^*A$ that is, 
$(mr)^*A\cong(qr)^*A$. We have $mr=\Ad(\d)$, $qr=1$ hence $\Ad(\d)^*A\cong A$ and 
$\Ad(\d\i)^*A\cong A$, as required.

Note also that for $(w,\l)$ as above we have:

(c) $\ff(\un{\dot{\Ad(\uD\i)^*\cl}}_{\e(w)}^\sh[\bd_w])\cong\ff(\udcl_w^\sh[\bd_w])$.
\nl
Indeed, let $K=\ff(\udcl_w^\sh[\bd_w])$. Clearly we have $m^*K\cong q^*K$ with $m,q$ as 
above. Then as in the proof of (b) we see that $\Ad(\d)^*K\cong K$. From the definitions 
we see that $\ff(\un{\dot{\Ad(\uD\i)^*\cl}}_{\e(w)}^\sh[\bd_w])=\Ad(\d\i)^*K$. Since
$\Ad(\d\i)^*K\cong K$, (c) follows.

\subhead 41.7\endsubhead
In this and next subsection we assume that $\kk$ is an algebraic closure of a finite field.
From 41.1(c) we see that $\x_D^*:\cd(Z_{\em,D})@>>>\cd(C)$ restricts to a functor
$\cd^{cs}(Z_{\em,D})@>>>\cd^{cs}(C)$ hence, as in 36.8, the $\ca$-linear map 
$gr(\x_D^*):\fK(Z_{\em,D})@>>>\fK(C)$ is well defined; from 41.1(c) we see also that 

(a) $gr(\x_D^*)(\udcl^\sh_w[\bd_w])=(-v)^\rr[w;\uD\l]$ 
\nl
for $w\in\WW,\l\in\ufs$ such that $w\uD\l=\l$ and $\cl\in\l$. From (a) we see that
$gr(\x_D^*)$ is injective with image equal to $\fK(C)^D$, the $\ca$-submodule of $\fK(C)$ 
spanned by $\{[w;\uD\l];w\in\WW,\l\in\ufs,w\uD\l=\l\}$ or equivalently by 
$\{[w;\uD\l]';w\in\WW,\l\in\ufs,w\uD\l=\l\}$. Thus, $gr(\x_D^*)$ defines an isomorphism 
$\et':\fK(Z_{\em,D})@>\si>>\fK(C)^D$. Let $\et=\et'{}\i$.

Let $n\in\NN^*_\kk$. Let $\fK(C)^D_n$ be the $\ca$-submodule of $\fK(C)$ 
spanned by $\{[w;\uD\l];w\in\WW,\l\in\ufs_n,w\uD\l=\l\}$ or equivalently by 
$\{[w;\uD\l]';w\in\WW,\l\in\ufs_n,w\uD\l=\l\}$. 

\mpb

Let $u,w\in\WW,\l\in\ufs_n$ be such that $w\uD\l=\l$ and let $\cl\in\l$. From 37.3(c) we 
see that the $\ca$-linear map $gr(\Ph_u):\fK(Z_{\em,D})@>>>\fK(Z_{\em,D})$ is well defined;
we denote it again by $\Ph_u$. From 40.10(a), 41.2(a) we have
$$[u\i;\l]'*[w;\uD\l]'{}^\sh*[\e_D(u);\uD(u\i\l)]'
=(v^2-1)^{2\rr}\et'\Ph_u\et([w;\uD\l]'{}^\sh),$$
equality in $\fK(C)$. 
If $\l'\in\ufs_n,\l'\ne\l$ we have (from 40.7) that
$[u\i;\l']'*[w;\uD\l]'{}^\sh*[\e_D(u);\uD(u\i\l')]'=0$. It follows that
$$(v^2-1)^{2\rr}\et'\Ph_u\et([w,\uD\l]'{}^\sh))
=\sum_{\l'\in\ufs_n}[u\i;\l']'*[w;\uD\l]'{}^\sh*[\e_D(u);\uD(u\i\l')]'$$
Using this and the definition of $\fK(C)^D_n$ we see that
$$(v^2-1)^{2\rr}\et'\Ph_u\et(x)=\sum_{\l'\in\ufs_n}[u\i;\l']'*x*[\e_D(u);\uD(u\i\l')]'$$
for any $x\in\fK(C)^D_n$. Applying $\et$ to both sides we obtain
$$(v^2-1)^{2\rr}\Ph_u\et'(x)=\sum_{\l'\in\ufs_n}\et([u\i;\l']'*x*[\e_D(u);\uD(u\i\l')]')
\tag b$$
for any $x\in\fK(C)^D_n$.

\subhead 41.8\endsubhead
In the setup of 41.4, let $A$ be a character sheaf on $Z_{J,D}$. From 36.9(b) we see that 
the condition that, if $(w',\uD\l')\in\WW\T\ufs$ is such that $w'\uD\l'=\l'$, then we have 
$A\dsv\ff(\udcl'_{w'}{}^\sh[\bd_{w'}])$ if and only if $A$ appears with coefficient $\ne0$ 
in the expansion of $\ff(\udcl'_{w'}{}^\sh[\bd_{w'}])\in\fK(Z_{J,D})$ as a linear 
combination of the canonical basis of $\fK(Z_{J,D})$. Hence from 41.5(a) we deduce:

(a) {\it There exists $(w,\uD\l)\in\boc_A$ such that $w\uD\l=\l$ and $A$ appears with 
non-zero coefficient in $\ff(\udcl_w^\sh[\bd_w])\in\fK(Z_{J,D})$. If 
$(w',\uD\l')\in\WW\T\ufs$ is such that $w'\uD\l'=\l'$ and $A$ appears with non-zero 
coefficient in $\ff(\udcl'_{w'}{}^\sh[\bd_{w'}])\in\fK(Z_{J,D})$ then 
$(w,\uD\l)\prq_{J,J'}(w',\uD\l')$. Here $\cl\in\l,\cl'\in\l'$.}
\nl
Clearly, property (a) characterizes $\boc_A$. 

\subhead 41.9\endsubhead
Let $J\sub J'\sub\II$ and let $D'$ be another connected component of $G$. Let 
$A_0\in\cd(Z_{J,D})$, $A'\in\cd(Z_{\e_D(J'),D'})$. We show:

(a) $\ff_{J,J'}(A_0)*A'\cong\ff_{J,J'}(A_0*\fe_{\e_D(J),\e_D(J')}A')$ in $\cd(Z_{J',D'D})$.
\nl
Indeed, from the definitions we see that both sides of (a) can be identified with
$b_!c^*(A_0\bxt A')$ where $b,c$ are as in the diagram

$Z_{J,D}\T Z_{\e_D(J'),D'}@<c<<Y@>b>>Z_{J',D'D}$
\nl
where 
$$\align Y=&\{(P,R,R',gU_R,g'U_{R'});
P\in\cp_J,R\in\cp_{J'},R'\in\cp_{\e_D(J')},\\&gU_R\in D/U_R,
g'U_{R'}\in D'/U_{R'},gRg\i=R',P\sub R\},\endalign$$

$c$ is $(P,R,R',gU_R,g'U_{R'})\m((P,gU_P),(R',g'U_{R'})$,

$b$ is $(P,R,R',gU_R,g'U_{R'})\m(R,g'gU_R)$.

An entirely similar proof shows that, if $A\in\cd(Z_{J',D})$, $A'_0\in\cd(Z_{\e_D(J),D'})$
then

(b) $A*\ff_{\e_D(J),\e_D(J')}(A'_0)\cong\ff_{J,J'}(\fe_{J,J'}A*A'_0)$ in $\cd(Z_{J',D'D})$.

\subhead 41.10\endsubhead
Let $\boc$ be a two-sided cell in $\WW\T\ufs$. Let $\bar\boc$ be the set of all 
$(w,\l)\in\WW\T\ufs$ such that $(w,\l)\prq_{\II,\II}(y,\nu)$ for some/any $(y,\nu)\in\boc$.

If $K\in\cd(Z_{\em,D})$, we say that $K\in\cd^{cs}_{\bar\boc}(Z_{\em,D})$ if for any 
$j\in\ZZ$ and simple subquotient $A$ of ${}^pH^j(K)$ satisfies $\boc_A\sub\bar\boc$.

Let $D'$ be another connected component of $G$. We show:

(a) {\it If $K\in\cd^{cs}_{\bar\boc}(Z_{\em,D})$, $K'\in\cd^{cs}(Z_{\e_D(J'),D'})$, then
$K*K'\in\cd^{cs}_{\bar\boc}(Z_{\em,D'D})$.}
\nl
We may assume that $\kk$ is an algebraic closure of a finite field. We may assume that 
$K\in\hZ_{\em,D}$ and $\boc_K\sub\bar\boc$. Then there exists $(w,\uD\l)\in\boc_K$ such 
that $w\uD\l=\l$, $K\dsv\ff(\udcl_w^\sh[\bd_w])$, $\cl\in\l$. It 
is enough to show that, if $\tA\in\hZ_{\em,D'D}$ is such that $\tA\dsv K*K'$ then 
$\boc_{\tA}\sub\bar\boc$. Since $\ff(\udcl_w^\sh[\bd_w])$ is a semisimple complex (see the
line after 41.4(e)) we have $\ff(\udcl_w^\sh[\bd_w])\cong K[m]\op\tK$ for some $m\in\ZZ$,
$\tK\in\cd(Z_{\em,D'D})$. It follows that 
$\ff(\udcl_w^\sh[\bd_w])*K'\cong K*K'[m]\op\tK*K'$ hence 
$\tA\dsv\ff(\udcl_w^\sh[\bd_w])*K'$. By 41.9(a) we have
$\ff(\udcl_w^\sh[\bd_w])*K'\cong\ff(\udcl_w^\sh[\bd_w]*\fe(K'))$ hence 
$\tA\dsv\ff(\udcl_w^\sh[\bd_w]*\fe(K'))$. We deduce that there exists $K'_0\in\hZ_{\em,D'}$
such that $\tA\dsv\ff(\udcl_w^\sh[\bd_w]*K'_0)$ and $K''_0\in\hZ_{\em,D'D}$ such that 
$K''_0\dsv\udcl_w^\sh[\bd_w]*K'_0$, $\tA\dsv\ff(K''_0)$. We then have
$\x_{D'D}^*K''_0[\rr]\dsv\x_{D'D}^*(\udcl_w^\sh[\bd_w]*K'_0)[\rr]$, hence, using 41.1(d),
$\x_{D'D}^*K''_0[\rr]\dsv(\x_D^*\udcl_w^\sh[\bd_w])*\x_{D'}K'_0)$. Setting 
$gr(\x_{D'D}^*K''_0[\rr])=[w_1,\uD'\uD\l_1]\in\fK(C)$ with $(w_1,\l_1)\in\WW\ufs$ we see, 
using 41.3(a) that $(w_1,\uD'\uD\l_1)\prq_{\II,\II}(w,\uD\l)$. From $\tA\dsv\ff(K''_0)$ we
see using 41.5(a) that $\boc_{\tA}\prq_{\II,\II}(w_1,\uD'\uD\l_1)$ (that is, some/any
element of $\boc_{\tA}$ is $\prq_{\II,\II}(w_1,\uD'\uD\l_1)$). Using the transitivity of
$\prq_{\II,\II}$ we see that $\boc_{\tA}\prq_{\II,\II}(w,\uD\l)$. This proves (a).

An entirely similar argument shows:

(b) {\it If $K\in\cd^{cs}(Z_{\em,D})$, $K'\in\cd^{cs}_{\bar\boc}(Z_{\e_D(J'),D'})$, then
$K*K'\in\cd^{cs}_{\bar\boc}(Z_{\em,D'D})$.}

\head 42. Duality and the functor $\ff_{\em,\II}$\endhead
\subhead 42.1\endsubhead
In this section we fix a connected component $D$ of $G$. We write $\e$ instead of
$\e_D:\WW@>>>\WW$. We write $\ff$ instead of $\ff_{\em,J}:\cd(Z_{\em,D}@>>>\cd(Z_{J,D})$. 
We assume that $\kk$ is an algebraic closure of a finite field.

Let $J\sub\II$ be such that $\e(J)=J$. Recall from 30.3 that
$V_{J,D}=\{(P,gU_P);P\in\cp_J,gU_P\in N_DP/U_P\}$. As in 30.4 (with $J'=\II$) we consider 
the diagram $V_{J,D}@<c<<V_{J,\II,D}@>d>>D$ where 
$V_{J,\II,D}=\{(P,g);P\in\cp_J,g\in N_DP\}$, $c$ is $(P,g)\m(P,gU_P)$ and $d$ is 
$(P,g)\m g$. Define $\tf_J:\cd(V_{J,D})@>>>\cd(D)$, $\te_J:\cd(D)@>>>\cd(V_{J,D})$ by 
$\tf_JA=d_!c^*A,\te_JA'=c_!d^*A'$. (In the notation of 30.4 we have $\tf_J=\tf_{J,\II}$,
$\te_J=\te_{J,\II}$.) Define $f_J:\cd(V_{J,D})@>>>\cd(D)$, $e_J:\cd(D)@>>>\cd(V_{J,D})$ by
$f_JA=\tf_JA[[\a_J/2]]$, $e_JA=\te_JA[[\a_J/2]]$ where $\a_J=\dim\cp_J$. (In the notation
of 30.4 we have $f_JA=f_{J,\II}A(\a_J/2)$, $e_JA=e_{J,\II}A(-\a_J/2)$. Thus, $f_J,e_J$ are
the same, up to a twist, as $f_{J,\II},e_{J,\II}$.) 

From 30.5 (with $J'=\II$) we see that for $A\in\cd(V_{J,D})$, $A'\in\cd(D)$ we have 
canonically

(a) $\Hom_{\cd(V_{J,D})}(e_JA',A)=\Hom_{\cd(D)}(A',f_JA)$.
\nl
Let $CS(V_{J,D}),CS(D)$ be as in 38.1. From 38.2, 38.3 we see that

(b) {\it $f_J,e_J$ restrict to functors $CS(V_{J,D})@>>>CS(D)$, $CS(D)@>>>CS(V_{J,D})$ 
denoted again by $f_J,e_J$.}
\nl
We show:

(c) {\it if $A\in CS(V_{J,D})$ comes from a pure complex of weight $0$ with respect to a
rational structure over a finite subfield of $\kk$ then $f_JA$ (naturally regarded as a
mixed complex) is pure of weight $0$.}
\nl
Indeed, the functor $c^*$ preserves pure complexes of weight $0$ since $c$ is smooth with 
connected fibres; the functor $d_!$ preserves pure complexes of weight $0$ since $d$ is 
proper (see \cite{\DE, 6.2.6}) and $[[\a_J]]$ also preserves pure complexes of weight $0$.

We show:

(d) {\it if $A'\in CS(D)$ comes from a pure complex of weight $0$ with respect to a 
rational structure over a finite subfield of $\kk$ then $e_JA'$ (naturally regarded as a
mixed complex) is pure of weight $0$.}
\nl
Using (b), it is enough to show that for any simple $A$ as in (c), the natural action of 
Frobenius on the vector space $\Hom_{\cd(V_{J,D})}(e_JA',A)$ has weight $0$. Using (a) we
see that it is enough to show that the natural action of Frobenius on the vector space 
$\Hom_{\cd(D)}(A',f_JA)$ has weight $0$. This follows from (c) using (b).

Define an imbedding $s:V_{J,D}@>>>Z_{J,D}$ by $(P,gU_P)\m(P,P,gU_P)$. From the definitions
we see that

(e) {\it $\tf_J:\cd(V_{J,D})@>>>\cd(D)$ is the composition
                                 $\cd(V_{J,D})@>s_!>>\cd(Z_{J,D})@>\ff_{J,\II}>>\cd(D)$,}

(f) {\it $\te_J:\cd(D)@>>>\cd(V_{J,D})$ is the composition
                                 $\cd(D)@>\fe_{J,\II}>>\cd(Z_{J,D})@>s^*>>\cd(V_{J,D})$.}
\nl
Let $Y=\{(B,B',gU_B)\in Z_{\em,D};\po(B,B')\in\WW_J\}$ and let $r:Y@>>>Z_{\em,D}$ be the
inclusion. From the definitions we have

(g) $s_!s^*\ff_{\em,J}=\ff_{\em,J}r_!r^*:\cd(Z_{\em,D})@>>>\cd(Z_{J,D})$.
\nl
Note that $V_{J,D}={}^1Z_{J,D}$, see 36.2; hence the "character sheaves" on 
$V_{J,D}={}^1Z_{J,D}$ are defined as in 36.8 and $\cd^{cs}(V_{J,D}=\cd^{cs}({}^1Z_{J,D})$ 
is defined as 36.8. In particular, $\fK(V_{J,D})=\fK({}^1Z_{J,D})$ is defined. Let 
$\fK_0(V_{J,D})=\op_A\ZZ A\sub\fK(V_{J,D})$ where $A$ runs through the character sheaves on
$V_{J,D}$ (up to isomorphism).

From (b) we see that $\tf_J,\te_J$ restrict to functors $\cd^{cs}(V_{J,D})@>>>\cd^{cs}(D)$,
$\cd^{cs}(D)@>>>\cd^{cs}(V_{J,D})$ hence the $\ca$-linear maps
$gr(\tf_J):\fK(V_{J,D})@>>>\fK(D)$, $gr(\te_J):\fK(D)@>>>\fK(V_{J,D})$ are well defined; we
denote them by $\tf_J,\te_J$. Define $f_J:\fK(V_{J,D})@>>>\fK(D)$ by 
$f_J=v^{-\a_J}\tf_J$ and $e_J:\fK(D)@>>>\fK(V_{J,D})$ by $e_J=v^{-\a_J}\te_J$. We show:

(h) {\it $f_J:\fK(V_{J,D})@>>>\fK(D)$, $e_J:\fK(D)@>>>\fK(V_{J,D})$ restrict to group
homomorphisms $\fK_0(V_{J,D})@>>>\fK_0(D)$, $\fK_0(D)@>>>\fK_0(V_{J,D})$ denoted again by
$f_J,e_J$.}
\nl
It is enough to prove the following statement. If $x$ is a canonical basis element of 
$\fK(V_{J,D})$ (resp. $\fK(D)$) then $f_J(x)$ (resp. $e_J(x)$) is an $\NN$-linear 
combination of canonical basis elements of $\fK(D)$ (resp. $\fK(V_{J,D})$). This is
immediate from (c), (d).

Now, one checks easily that $r_!r^*:\cd(Z_{\em,D})@>>>\cd(Z_{\em,D})$ restricts to a 
functor $\cd^{cs}(Z_{\em,D})@>>>\cd^{cs}(Z_{\em,D})$. (Note that, if $w\in\WW,\l\in\ufs$,
$\cl\in\l$ and $w\uD\l=\l$, then $r_!r^*(\udcl_w)=\udcl_w$ for $w\in\WW_J$ and 
$r_!r^*(\udcl_w)=0$ for $w\in\WW-\WW_J$.) It follows that the $\ca$-linear map 
$gr(r_!r^*):\fK(Z_{\em,D})@>>>\fK(Z_{\em,D})$ (denoted by $\r_J$) is well defined. 

Let $\fK(C)^D,\et$ be as in 41.7. Define an $\ca$-linear map $\tir_J:\fK(C)^D@>>>\fK(C)^D$
by $[w;\uD\l]'\m[w;\uD\l]'$ if $w\in\WW_J,\l\in\ufs,w\uD\l=\l$ and $[w;\uD\l]'\m0$ if 
$w\in\WW-\WW_J,\l\in\ufs,w\uD\l=\l$. From the definitions we see that

(i) $\r_J\et(x)=\et\tir_J(x)$ for all $x\in\fK(C)^D$.

\subhead 42.2\endsubhead
We define an $\ca$-linear map $\dd:\fK(D)@>>>\fK(D)$ by
$$\dd(x)=\sum_{J;J\sub\II;\e(J)=J}(-1)^{|J_\e|}f_Je_J(x)$$
where $f_J,e_J$ are as in 42.1(h) and $J_\e$ is as in 38.1. We show:

(a) {\it Let $A$ be a character sheaf on $D$. Then $\dd(A)=\pm A'$ where $A'$ is a 
character sheaf on $D$. Moreover $\pm$ and $A'$ are the same as in 38.11(a).}
\nl
For any $J\sub\II$ such that $\e(J)=J$ let $\ck(V_{J.D})$ be as in 38.9. We shall identify
$\fK(V_{J,D})/(v-1)\fK(V_{J,D})=\ck(V_{J,D})$ as abelian groups in such a way that the
image of $A_1$ (a character sheaf on $V_{J,D}$) in $\fK(V_{J,D})/(v-1)\fK(V_{J,D})$ is
identified with the basis element $A_1$ of $\ck(V_{J,D})$. From the definitions we see that
the homomorphisms 

$\fK(D)/(v-1)\fK(D)@>>>\fK(V_{J,D})/(v-1)\fK(V_{J,D})@>>>\fK(D)/(v-1)\fK(D)$
\nl
induced by $e_J,f_J$ in 42.1(h) are then identified with the homomorphisms

$e_{J,\II}:\ck(D)@>>>\ck(V_{J,D})$,  $f_{J,\II}:\ck(V_{J,D})@>>>\ck(D)$ 
\nl
in 38.2, 38.3. It follows that the 
endomorphism of $\fK(D)/(v-1)\fK(D)$ induced by $\dd:\fK(D)@>>>\fK(D)$
is identified with the homomorphism $\ck(D)@>>>\ck(D)$ denoted in 38.10(a), 38.11 again by
$\dd$. Hence we have $\dd(A)=\pm A'+(v-1)x$ (in $\fK(D)$) 
where $\pm,A'$ are as in 38.11(a) 
and $x\in\fK(D)$. From 42.1(h) we see that $\dd(A)\in\fK_0(D)$. Since $\pm A'\in\fK_0(D)$,
we see that $(v-1)x\in\fK_0(D)$. Since $\fK_0(D)\cap(v-1)\fK(D)=0$, we have $(v-1)x=0$
and $x=0$. This proves (a).

\subhead 42.3\endsubhead
We have $H=H_D\op H'_D$ where $H_D$ (resp. $H'_D$) 
is the $\ca$-submodule of $H_n$ spanned by 
$\{\tT_w1_{\uD\l};w\in\WW,\l\in\ufs,w\uD\l=\l\}$ (resp. by 
$\{\tT_w1_{\uD\l};w\in\WW,\l\in\ufs,w\uD\l\ne\l\}$). Equivalently,

$H_D=\sum_{\l\in\ufs}1_\l H1_{\uD\l}\sub H$,
$H'_D=\sum_{\l,\l'\in\ufs;\l\ne\l'}1_{\l'}H1_{\uD\l}\sub H$.

Recall that $\o:\fK(C)@>\si>>H$ is defined in 40.7(b). Define an $\ca$-linear map 
$\tio:H@>>>\fK(C)^D$ by

$\tio(y)=\o\i(y)$ if $y\in H_D$,,

$\tio(y)=0$ if $y\in H'_D$.
\nl
Then $\et\tio(y)\in\fK(Z_{\em,D})$ is well defined for any $y\in H$. Here $\et$ is as 
in 41.7. 

Let $n\in\NN^*_\kk$. Let $H_{n,D}=H_D\cap H_n$. Note that $H_{n,D}$ is the $\ca$-submodule
of $H_n$ spanned by $\{\tT_w1_{\uD\l};w\in\WW,\l\in\ufs_n,w\uD\l=\l\}$.

For $J\sub\II$ such that $\e(J)=J$ we define an $\ca$-linear map
$\r_{J,n}:H_{n,D}@>>>H_{n,D}$ by 

$\tT_w1_{\uD\l}\m\tT_w1_{\uD\l}$ if $w\in\WW_J,\l\in\ufs_n,w\uD\l=\l$,

$\tT_w1_{\uD\l}\m0$ if $w\in\WW-\WW_J,\l\in\ufs_n,w\uD\l=\l$.
\nl
We have the following result.

\proclaim{Lemma 42.4} For any $y\in H_{n,D}$ we have $\dd(\ff\et\tio(y))=\ff\et\tio(\d(y))$
where $\d=\sum_{J\sub\II;\e(J)=J}(-1)^{|J_\e|}\d_J$ with $\d_J:H_{n,D}@>>>H_{n,D}$ given by
\lb $\d_J(y)=\r_{J,n}(\sum_{u\in\WW^J}\tT_{u\i}y\tT_{\e_D(u)})$ (the sum in the right hand
side is computed in $H_n$ but it belongs to $H_{n,D}$).
\endproclaim
Applying 37.2 with $K,K',J$ repaced by $\em,J,\II$ and with $A'\in\cd^{cs}(Z_{\em,D})$ we 
obtain 

$\fe_{J,\II}\ff A'\Bpq\{\ff_{\em,J}\Ph_uA'[[-m_u]];u\in\WW^J\}$ 
\nl
(in $\cd(Z_{J,D})$, with $\Ph_u:\cd(Z_{\em,D})@>>>\cd(Z_{\em,D})$ as in 37.1 and 
$m_u=\a_J-\l(u)$ where $\a_J=\dim\cp_J$. Applying here $s^*$ we obtain

$s^*\fe_{J,\II}\ff A'\Bpq\{s^*\ff_{\em,J}\Ph_uA'[[-m_u]];u\in\WW^J\}$.
\nl
We replace $s^*\fe_{J,\II}$ by $\te_J$ (see 42.1(f)) and we apply $\tf_J=\ff_{J,\II}s_!$
(see 42.1(e)); we obtain

$\tf_J\te_J\ff A'\Bpq\{\ff_{J,\II}s_!s^*\ff_{\em,J}\Ph_uA'[[-m_u]];u\in\WW^J\}$.
\nl
Using now 42.1(g) we obtain

$\tf_J\te_J\ff A'\Bpq\{\ff_{J,\II}\ff_{\em,J}r_!r^*\Ph_uA'[[-m_u]];u\in\WW^J\}$.
\nl
Here we replace $\ff_{J,\II}\ff_{\em,J}$ by $\ff$ (see 36.4(b)). This (or rather 
its mixed analogue) gives rise to the following equality in $\fK(D)$:
$$\tf_J\te_J\ff(x')=\sum_{u\in\WW^J}v^{2m_u}\ff\r_J\Ph_u(x')$$
for any $x'\in\fK(Z_{\em,D})$, or equivalently
$$f_Je_J\ff(x')=\sum_{u\in\WW^J}v^{2m_u-2\a_J}\ff\r_J\Ph_u(x').$$
Taking $x'=\et(x)$ where $x\in\fK(C)^D_n$ (see 41.7) and using 41.7(b) we obtain
$$\align&(v^2-1)^{2\rr}f_Je_J\ff\et(x)\\&=\sum_{u\in\WW^J}\sum_{\l\in\ufs_n}
v^{-2l(u)}\ff\r_J\et([u\i;\l]'*x*[\e_D(u);\uD(u\i\l)]')\endalign$$
and using 42.1(i),
$$\align&(v^2-1)^{2\rr}f_Je_J\ff\et(x)\\&=\sum_{u\in\WW^J}\sum_{\l\in\ufs_n}
v^{-2l(u)}\ff\et\tir_J([u\i;\l]'*x*[\e_D(u);\uD(u\i\l)]')\endalign$$
for any $x\in\fK(C)^D_n$. Here we replace $x$ by $\tio(y)$ where $y\in H_{n,D}$ and 
$\tir_J|_{\fK(C)^D_n}$ by $\tio|_{H_{n,D}}\r_{J,n}\o_{\fK(C)^D_n}$; using 40.7(b) we 
obtain:
$$\align&f_Je_J\ff\et\tio(y)=\sum_{u\in\WW^J}\sum_{\l\in\ufs_n}
\ff\et\tio\r_{J,n}(\tT_{u\i}1_\l y\tT_{\e_D(u)}1_{\uD(u\i\l)})\\&=
\ff\et\tio\r_{J,n}(\sum_{u\in\WW^J}\tT_{u\i}y\tT_{\e_D(u)}).\endalign$$
The lemma is proved.

\subhead 42.5\endsubhead
As in 34.12 let $\fU$ be the subfield of $\bbq$ generated by the roots of $1$. Let 
$\Ph:H_n^D@>>>\ca\ot_\ZZ H_n^{D,\iy}$ be as in 34.12 (a special case of a definition in 
34.1) and let $\Ph^1:H_n^{D,1}@>>>\fU\ot_\ZZ H_n^{D,\iy}$ be the specialization of $\Ph$ 
for $v=1$ (see 34.12(b)). Let $\tca=\fU[v,v\i]$, let $H_n^{D,\tca}=\tca\ot\ca H_n^D$ and 
let $\Ph^{\tca}:H_n^{D,\tca}@>>>\tca\ot_\ZZ H_n^{D,\iy}$ be the homomorphism obtained from 
$\Ph$ by extending the scalars from $\ca$ to $\tca$.

Let $E$ be an $H_n^{D,1}$-module of finite dimension over $\fU$. Since $\Ph^1$ is an 
isomorphism of $\fU$-algebras (see 34.12(b)) we may regard $E$ as an 
$\fU\ot_\ZZ H_n^{D,\iy}$-module $E^\iy$ via $(\Ph^1)\i$. By extension of scalars,
$\tca\ot_\fU E^\iy$ is naturally a module over 

$\tca\ot_\fU(\fU\ot_\ZZ H_n^{D,iy})=\tca\ot_\ZZ H_n^{D,iy}$
\nl
and this can be regarded as an $H_n^{D,\tca}$-module $E^{\tca}$ via $\Ph^{\tca}$.

Let $J\sub\II$ be such that $\e(J)=J$. Let $H_{J,n}^D$ be the $\ca$-algebra of $H_n^D$
generated by $1_\l,\l\in\ufs_n$ by $\tT_w,w\in\WW_J$ and by $\tT_{\uD}$. Note that 
$\{\tT_{w\uD'}1_\l;w\in\WW_J,\uD'=\text{power of }\uD\}$ is an $\ca$-basis of $H_{J,n}^D$.
Let $H_{J,n}^{D,1}=\fU\ot_\ca H_{J,n}^D$ where $\fU$ is regarded as an $\ca$-algebra via 
$v\m1$. Let $H_{J,n}^{D,\tca}=\tca\ot_\ca H_{J,n}^D$. Note that $H_{J,n}^{D,\tca}$ is 
naturally a subalgebra of $H_n^{D,\tca}$. Hence $E^{\tca}$ may be regarded as an
$H_{J,n}^{D,\tca}$-module $(E^{\tca})_J$. This $H_{J,n}^{D,\tca}$-module may be induced to
an $H_n^{D,\tca}$-module 
$\IND((E^{\tca})_J):=H_n^{D,\tca}\ot_{H_{J,n}^{D,\tca}}E^{\tca}_J$.

Next, $H_{J,n}^{D,1}$ is naturally a subalgebra of $H_n^{D,1}$. Hence $E$ may be regarded 
as an $H_{J,n}^{D,1}$-module $E_J$. This $H_{J,n}^{D,1}$-module may be induced to an
$H_n^{D,1}$-module $\ind(E_J):=H_n^{D,1}\ot_{H_{J,n}^{D,1}}E_J$. Define an 
$H_{J,n}^{D,\tca}$-module $(\ind(E_J))^{\tca}$ in terms of $\ind(E_J)$ in the same way as 
$E^{\tca}$ was defined in terms of $E$. By extension of scalars from $\tca$ to $\fU(v)$ 
(the quotient field of $\tca$), $\IND((E^{\tca})_J)$, $(\ind(E_J))^{\tca}$ give rise to
$\fU(v)\ot_\ca H_n^D$-modules $\fU(v)\ot_{\tca}\IND((E^{\tca})_J)$, 
$\fU(v)\ot_{\tca}(\ind(E_J))^{\tca}$. We show:

(a) {\it These two $\fU(v)\ot_\ca H_n^D$-modules are isomorphic.}
\nl
Since $\fU(v)\ot_\ca H_n^D$, $H_n^{D,1}$ are (finite dimensional) semisimple algebras (see 
34.12) it follows by standard arguments that it is enough to show that 
$\IND((E^{\tca})_J)$,
$(\ind(E_J))^{\tca}$ become isomorphic $H_n^{D,1}$-modules under the specialization $v=1$.
First we note that under the specialization $v=1$, $E^{\tca}$ becomes the 
$H_n^{D,1}$-module $E$. (This is because the specialization of $\Ph^{\tca}$ at $v=1$ 
cancels $(\Ph_1)\i$.) In particular, the specialization of $(\ind(E_J))^{\tca}$ for $v=1$ 
is $\ind(E_J)$. Moreover, from the definition of induction, the specialization of 
$\IND((E^{\tca})_J)$ for $v=1$ is the same as $\ind(E'_J)$ where $E'$ is the specialization
of $E^{\tca}$ for $v=1$ that is, $E'=E$. This proves (a).

\proclaim{Lemma 42.6}We preserve the setup of 42.5. Let $\fU(v)\ot_{\tca}E^{\tca}$, 
$\fU(v)\ot_{\tca}(\ind(E_J))^{\tca}$ be the $\fU(v)\ot_\ca H_n^D$-module obtained from 
$E^{\tca}$, $(\ind(E_J))^{\tca}$ by extension of scalars from $\tca$ to $\fU(v)$. Let 
$y\in H_{n,D}$. We have:

$\tr(\d_J(y)\tT_{\uD},\fU(v)\ot_{\tca}E^{\tca})
=\tr(y\tT_{\uD},\fU(v)\ot_{\tca}(\ind(E_J))^{\tca})$.
\endproclaim
Let $H_{J,n}$ be the $\ca$-subalgebra of $H_n$ defined in 31.8. Define an $\ca$-linear map
$p_J:H_n@>>>H_{J,n}$ by $p_J(\tT_z1_\l))=\tT_z1_\l$ if $z\in\WW_J,\l\in\ufs_n$, 
$p_J(\tT_z1_\l)=0$ if $z\in\WW-\WW_J,\l\in\ufs_n$. We show that

(a) $p_J(\tT_uh')=\d_{u,1}h'$ if $u\in\WW^J,h'\in H_{J,n}$.
\nl
We may assume that $h'=\tT_b1_\l,b\in\WW_J,\l\in\ufs_n$. Then
$p_J(\tT_u\tT_b1_\l)=p_J(\tT_{ub}1_\l)=\d_{u,1}\tT_{ub}1_\l=\d_{u,1}\tT_b1_\l$, as 
required. 

We show:

(b) $p_J(hh')=p_J(h)h'$ for any $h\in H_n,h'\in H_{J,n}$.
\nl
We may assume $h=\tT_u\tT_b1_\nu,h'=\tT_a1_\l$, $u\in\WW^J,a,b\in\WW_J$, $\l,\nu\in\ufs_n$.
We must show tha $p_J(\tT_u\tT_b1_\nu\tT_a1_\l)=p_J(\tT_u\tT_b1_\nu)\tT_a1_\l$. If $u\ne1$,
both sides are zero by (a). If $u=1$ both sides are $\tT_b1_\nu\tT_a1_\l$. This proves (b).

By 34.13(a) we have 

(c) $p_\em(\tT_x\tT_{x'}1_\l)=\d_{xx',1}$ for $x,x'\in\WW,\l\in\ufs_n$.
\nl
For $u,u'\in\WW^J,\l\in\ufs_n$ we write $\tT_{u\i}\tT_{u'}1_\l=\sum_{a\in\WW}f_a\tT_a1_\l$
where $f_a\in\ca$. For $a'\in\WW_J$ we have 
$$\tT_{a'{}\i u\i}\tT_{u'}1_\l=\tT_{a'{}\i}\tT_{u\i}\tT_{u'}1_\l
=\sum_{a\in\WW}f_a\tT_{a'{}\i}\tT_a1_\l.$$
Applying $p_\em$ to this and using (c) gives $f_{a'}=\d_{u',ua'}=\d_{a',1}\d_{u,u'}$ so 
that
 
$p_J(\tT_{u\i}\tT_{u'}1_\l)=\sum_{a\in\WW_J}f_a\tT_a1_\l=\d_{u,u'}\tT_11_\l$.
\nl
Since this holds for any $\l\in\ufs_n$ we have

(d) $p_J(\tT_{u\i}\tT_{u'})=\d_{u,u'}\tT_1$.
\nl
Let $w\in\WW,\l\in\ufs_n,u\in\WW^J$. We have 

$\tT_w1_\l\tT_u=\sum_{u'\in\WW^J,a\in\WW_J}c_{w,u,u',a,\l}\tT_{u'}\tT_a1_{u\i\l}$
\nl
where $c_{w,u,u',a,\l}\in\ca$ are uniquely determined. It follows that
$$\tT_{u\i}\tT_w1_\l\tT_{\e(u)}
=\sum_{u'\in\WW^J,a\in\WW_J}c_{w,\e(u),u',a,\l}\tT_{u\i}\tT_{u'}\tT_a1_{\e(u)\i\l}.$$
Applying $p_J$ and using (b),(d) we obtain
$$\align&p_J(\tT_{u\i}\tT_w1_\l\tT_{\e(u)})
=\sum_{u'\in\WW^J,a\in\WW_J}c_{w,\e(u),u',a,\l}p_J(\tT_{u\i}\tT_{u'})\tT_a1_{\e(u)\i\l}\\&
=\sum_{u'\in\WW^J,a\in\WW_J}c_{w,\e(u),u',a,\l}\d_{u,u'}\tT_a1_{\e(u)\i\l}
=\sum_{a\in\WW_J}c_{w,\e(u),u,a,\l}\tT_a1_{\e(u)\i\l}.\tag e\endalign$$
Let $(e_i)_{i\in X}$ be a basis of the free $\tca$-module $E^{\tca}$. For 
$a\in\WW_J,\l\in\ufs_n$ we have 
$\tT_a1_\l\tT_{\uD}e_i=\sum_{i'\in X}\tc_{a,\l,i,i'}e_{i'}$ where $\tc_{a,\l,i,i'}\in\tca$.

Since $H_n^{D,\tca}$ is a free right $H_{J,n}^{D,\tca}$-module with basis 
$\{\tT_u;u\in\WW^J\}$ we see that $\{\tT_u\ot e_i;u\in\WW^J,i\in X\}$ is a basis of the 
free $\tca$-module $\ind((E^{\tca})_J)$.

Let $w\in\WW,\l\in\ufs_n,u\in\WW^J$ be such that $w\uD\l=\l$. In $\IND((E^{\tca})_J)$ we 
have
$$\align&\tT_w1_\l\tT_{\uD}(\tT_u\ot e_i)=(\tT_w1_\l\tT_{\e(u)}\tT_{\uD})\ot e_i\\&
=\sum_{u'\in\WW^J,a\in\WW_J}c_{w,\e(u),u',a,\l}
(\tT_{u'}\tT_a1_{\e(u)\i\l}\tT_{\uD})\ot e_i\\&
=\sum_{u'\in\WW^J,a\in\WW_J}c_{w,\e(u),u',a,\l}\tT_{u'}\ot(\tT_a1_{\e(u)\i\l}\tT_{\uD}e_i)
\\&=\sum_{u'\in\WW^J,a\in\WW_J,i'\in X}
c_{w,\e(u),u',a,\l}\tc_{a,\e(u)\i\l,i,i'}\tT_{u'}\ot e_{i'}.\endalign$$
Hence, using (e),
$$\align&\tr(\tT_w1_\l\tT_{\uD},\IND((E^{\tca})_J))
=\sum_{u\in\WW^J,a\in\WW_J,i\in X}c_{w,\e(u),u,a,\l}\tc_{a,\e(u)\i\l,i,i}\\&
=\sum_{u\in\WW^J,a\in\WW_J}c_{w,\e(u),u,a,\l}\tr(\tT_a1_{\e(u)\i\l}\tT_{\uD},E^{\tca})\\&
=\sum_{u\in\WW^J}\tr(\sum_{a\in\WW_J}c_{w,\e(u),u,a,\l}\tT_a1_{\e(u)\i\l}
\tT_{\uD},E^{\tca})
\\&=\tr(p_J(\sum_{u\in\WW^J}\tT_{u\i}\tT_w1_\l\tT_{\e(u)})\tT_{\uD},E^{\tca})\\&
=\tr(\r_{J,n}(\sum_{u\in\WW^J}\tT_{u\i}\tT_w1_\l\tT_{\e(u)})\tT_{\uD},E^{\tca})
=\tr(\d_J(\tT_w1_\l)\tT_{\uD},E^{\tca}).\endalign$$
Thus we have

$\tr(\d_J(\tT_w1_\l)\tT_{\uD},E^{\tca})=\tr(\tT_w1_\l\tT_{\uD},\IND((E^{\tca})_J))=
\tr(\tT_w1_\l\tT_{\uD},(\ind(E_J))^{\tca})$
\nl
where the second equality follows from 42.5(a). Since the elements $\tT_w1_\l$ as above 
generate the $\ca$-module $H_{n,D}$, the lemma follows.

\subhead 42.7\endsubhead
Let $\cv$ be the $\QQ$-vector subspace of $\QQ\ot\Hom(\kk^*,\TT)$ spanned by the coroots. 
Let $\cv_\RR=\RR\ot_\QQ\cv$. The kernels of the roots $\cv_\RR@>>>\RR$ a hyperplane 
arrangement which defines a partition of $\cv_\RR$ into facets in a standard way. Let $\cf$
be the set of facets. Now the orbits of $\WW$ on $\cf$ are naturally indexed by the various
subsets $J$ of $\II$. This gives a partition $\cf=\sqc_{J\sub\II}\cf_J$. For example 
$\cf_\em$ consists of all Weyl chambers. If $F\in\cf_J$ then $F$ is homeomorphic to a real
affine space of dimension $|\II-J|$ hence we have $H^i_c(F)=0$ if $i\ne|\II-J|$ and 
$H^{|\II-J|}_c(F)=\L^{|\II-J|}[F]$; here we write $H^i_c(?)$ instead of $H^i_c(?,\RR)$, 
$[F]$ denotes the vector subspace of $\cv_\RR$ in which $F$ is open dense and 
$\L^{|\II-J|}[F]$ is the top exterior power of $[F]$. Note that $[F]=\RR\ot_\QQ([F]_\QQ)$ 
for a well defined $\QQ$-subspace $[F]_\QQ$ of $\cv$. For any $\uD$-orbit $\co$ on the set
of subsets of $\II$ let $\cv_\RR^\co=\cup_{J\in\co}\cup_{F\in\cf_J}F\sub\cv_\RR$ and let 
$r_\co=|\II-J|$ for some/any $J\in\co$. We have $H^i_c(\cv_\RR^\co)=0$ if $i\ne r_\co$, 
$H^{r_\co}_c(\cv_\RR^\co)=\op_{J\in\co}\op_{F\in\cf_J}\L^{r_\co}[F]$. Note also that 
$H^i_c(\cv_\RR)=0$ if $i\ne|\II|$ and $H^{|\II|}_c(\cv_\RR)=\L^{|\II|}\cv_\RR$. The 
$\WW^D$-action on $\TT$ induces a linear action of $\WW^D$ on $\cv_\RR$. This action 
restricts for any $\co$ to a $\WW^D$-action on $\cv_\RR^\co$ and this induces a 
$\WW^D$-action on $H^{r_\co}_c(\cv_\RR^\co)$. It also induces a natural $\WW^D$-action on 
$H^{|\II|}_c(\cv_\RR)=\L^{|\II|}\cv_\RR$. The long cohomology exact sequences attached to 
the partition $\cv_\RR=\cup_\co\cv_\RR^\co$ show that
$(-1)^{|\II|}H^{|\II|}_c(\cv_\RR)=\sum_{\co}(-1)^{r_\co}H^{r_\co}_c(\cv_\RR^\co)$ in the 
Grothendieck group of representations of $\WW^D$ over $\RR$ that is,
$$\align&\L^{|\II|}\cv_\RR
\op\op_{\co;r_\co=|\II|+1\mod2}(\op_{J\in\co}\op_{F\in\cf_J}\L^{r_\co}[F])\\&
\cong\op_{\co;r_\co=|\II|\mod2}(\op_{J\in\co}\op_{F\in\cf_J}\L^{r_\co}[F])\endalign
$$
as representations of $\WW^D$ over $\RR$. All real representations in this formula come 
naturally from representations of $\WW^D$ over $\QQ$. Hence the previous formula remains 
valid (as representations of $\WW^D$ over $\QQ$) if $\cv_\RR,[F]$ are replaced by 
$\cv,[F]_\QQ$ and the exterior powers are taken over $\QQ$. Tensoring both sides (over 
$\QQ$) by $\fU$ (as in 42.5) we obtain
$$\align&\L^{|\II|}\cv_\fU\op
\op_{\co;r_\co=|\II|+1\mod2}(\op_{J\in\co}\op_{F\in\cf_J}\L^{r_\co}[F]_\fU)\\&
\cong\op_{\co;r_\co=|\II|\mod2}(\op_{J\in\co}\op_{F\in\cf_J}\L^{r_\co}[F]_\fU)\tag a
\endalign$$
as representations of $\WW^D$ over $\fU$; here $\cv_\fU=\fU\ot_\QQ\cv$,
$[F]_\fU=\fU\ot_\QQ[F]_\QQ$ and the exterior powers are taken over $\fU$. We may view (a) 
as an isomorphism of $H_n^{D,1}$-modules: the $\WW^D$-modules in (a) may be viewed as 
$H_n^{D,1}$-modules via the algebra homomorphism $H_n^{D,1}@>>>\fU[\WW^D]$ given by 
$\tT_w\m w$ for $w\in\WW^D$, $1_\l\m0$ for $\l\ne\l_0$, $1_{\l_0}\m1$ (here $\l_0$ is the 
neutral element of the abelian group $\ufs_n$, see 28.1). 

We define an $\fU$-linear map $\D:H_n^{D,1}@>>>H_n^{D,1}\ot H_n^{D,1}$ by
$\D(\tT_w)=\tT_w\ot\tT_w$ for $w\in\WW^D$ and
$\D(1_\l)=\sum_{\l_1,\l_2\in\ufs_n;\l_1\l_2=\l}1_{\l_1}\ot1_{\l_2}$ for $\l\in\ufs_n$.
(Here we use the abelian group structure on $\ufs_n$, a subgroup of $\ufs$, see 28.1.) This
makes $H_n^{D,1}$ into a Hopf algebra. (Note that the analogous formulas do not make 
$H_n^D$ into a Hopf algebra.) It follows that for any two $H_n^{D,1}$-modules $E_1,E_2$, 
the $\fU$-vector space $E_1\ot E_2$ is naturally an $H_n^{D,1}$-module. 

Now let $E$ be an $H_n^{D,1}$-module of finite dimension over $\fU$. Then we can take 
tensor product of each $H_n^{D,1}$-module in (a) with $E$ and we obtain an isomorphism of 
$H_n^{D,1}$-modules
$$E\ot\L^{|\II|}\cv_\fU\op 
\op_{\co;r_\co=|\II|+1\mod2}X_\co\cong\op_{\co;r_\co=|\II|\mod2}X_\co$$
where $X_\co=E\ot\op_{J\in\co}\op_{F\in\cf_J}\L^{r_\co}[F]_\fU$. Applying to this the 
functor $E\m E^{\tca}$, see 42.5, we deduce an isomorphism of $H_n^{D,\tca}$-modules
$$(E\ot\L^{|\II|}\cv_\fU)^{\tca}\op 
\op_{\co;r_\co=|\II|+1\mod2}X_\co^{\tca}\cong\op_{\co;r_\co=|\II|\mod2}X_\co^{\tca}.$$

We deduce that for $y\in H_{n,D}$ we have
$$\tr(y\tT_{\uD},(E\ot\L^{|\II|}\cv_\fU)^{\tca})=
\sum_{\co}(-1)^{r_\co+|\II|}\tr(y\tT_{\uD},X_\co^{\tca}).\tag b$$
We have $X_\co=\op_{J\in\co}X^J$ where $X^J=E\ot(\op_{F\in\cf_J}\L^{r_\co}[F]_\fU)$. 

Assume first that $\co$ consists of at least two subsets of $\II$. Then each $X_J$ is 
stable under $H_n^{D,1}$ and is mapped by $\tT_{\uD}$ into $X_{J'}$ with $J\ne J'$. From 
the definitions we have $X_\co^{\tca}=\op_{J\in\co}\tca\ot_\fU X_J$ as an $\tca$-module and
each summand $\tca\ot_\fU X_J$ is stable under $H_n$ and is mapped by $\tT_{\uD}$ into a 
summand $\tca\ot_\fU X_{J'}$ with $J\ne J'$. It follows that for our $\co$ we have

(c) $\tr(y\tT_{\uD},X_\co^{\tca})=0$.
\nl
Next assume that $\co$ consists of a single subset $J$ of $\II$. We have $\uD(J)=J$. Let 
$F_J$ be the unique facet in $\cf_J$ such that $F_J$ is contained in the closure of the 
dominant Weyl chamber. Then $F_J$ is stable under the the subgroup $\WW_J^D$ of $\WW^D$
generated by $\WW_J$ and $\uD$ and $X_\co$ may be identified with 
$E\ot(H_n^{D,1}\ot_{H_{J,n}^{D,1}}(\L^{|\II-J|}[F_J]_\fU))$. Here $\L^{|\II-J|}[F_J]_\fU$ 
is regarded as a $WW_J^D$-module and then is viewed as a $H_{J,n}^{D,1}$-module via the 
canonical algebra homomorphism $H_{J,n}^{D,1}@>>>\fU[\WW_J^D]$; thus $1_\l$ acts on it as 
$1$ if $\l=\l_0$ and as $0$ if $\l\ne\l_0$. Note that in the $\WW_J^D$-module
$\L^{|\II-J|}[F_J]_\fU$, $\WW_J$ acts trivially (since $\WW_J$ acts trivially on 
$[F_J]_\fU$) and $\uD$ acts as multiplication by $(-1)^{|\II-J|-|(\II-J)_\e|}$). Let 
$X'_\co=E\ot(H_n^{D,1}\ot_{H_{J,n}^{D,1}}\fU)$ where $\fU$ is regarded as a
$H_{J,n}^{D,1}$-module coming from the trivial representation of $\WW_J^D$. We see that we
may identify $X_\co,X'_\co$ in a way compatible with the $H_n^1$-module structures and so
that the action of $\tT_{\uD}$ on $X_\co$ corresponds to $(-1)^{|\II-J|-|(\II-J)_\e|}$ 
times the action of $\tT_{\uD}$ on $X'_\co$. Using the definitions we see that we may 
identify $X_\co^{\tca},X'_\co{}^{\tca}$ in a way compatible with the $H_n$-module 
structures and so that the action of $\tT_{\uD}$ on $X_\co^{\tca}$ corresponds to 
$(-1)^{|\II-J|-|(\II-J)_\e|}$ times the action of $\tT_{\uD}$ on $X'_\co{}^{\tca}$. From
the definitions we have $X'_\co=\ind(E_J)$ (notation of 42.5). We see that for our $\co$ we
have
$$\tr(y\tT_{\uD},X_\co^{\tca})=
(-1)^{|\II-J|-|(\II-J)_\e|}\tr(y\tT_{\uD},(\ind(E_J))^{\tca}).\tag d$$
From the definitions (34.4) we see that there is a unique $\tca$-algebra homomorphism 
$\vt:H_n^{D,\tca}D@>>>H_n^{D,\tca}$ such that 

$\vt(1_\l)=1_\l$ for any $\l\in\ufs_n$,

$\vt(\tT_w)=(-1)^{l(w)}\tT_{w\i}\i$ for any $w\in\WW$,

$\vt(\tT_{\uD})=(-1)^{|\II|-|\II_\e|}\tT_{\uD}$.
\nl
We have $\vt^2=1$. 

Using $\vt$ and $E^{\tca}$ we can define a new $H_n^{D,\tca}$-module $E^{\tca,\vt}$ with 
the same underlying $\tca$-module as $E^{\tca}$ but with $x\in H_n^{D,\tca}$ acting on
$E^{\tca,\vt}$ in the same way that $\vt(x)$ acts on $E^{\tca}$. We show: 

(e) {\it under extension of scalars from $\tca$ to $\fU(v)$, the $H_n^{D,\tca}$-modules 
$E^{\tca,\vt}$ and $(E\ot\L^{|\II|}\cv_\fU)^{\tca}$ become isomorphic
$\fU(v)\ot_\ca H_n^D$-modules.}
\nl
As in the proof of 42.5(a) it is enough to show that these $H_n^{D,\tca}$-modules become
isomorphic $H_n^{D,1}$-modules under the specialization $v=1$. Thus it is enough to show 
that $E^{\tca,\vt}|_{v=1}\cong E\ot\L^{|\II|}\cv_\fU$ as $H_n^{D,1}$-modules. Now the
underlying $\fU$-vector space of $E^{\tca,\vt}|_{v=1}$ is $E$ but the action of
$x\in H_n^{D,1}$ on $E^{\tca,\vt}|_{v=1}$ is the same as the action of $\vt_1(x)$ on $E$. 
Here $\vt_1:H_n^{D,1}@>>>H_n^{D,1}$ is the specialization of $\vt_1$ for $v=1$. Note that 
$\vt_1(1_\l)=1_\l$ for any $\l\in\ufs_n$ and $\vt_1(\tT_w)=\g_w\tT_w$ for any $w\in\WW^D$, 
where $\g_w=\pm1$ is the scalar by which $w$ acts in the $\WW^D$-module 
$\L^{|\II|}\cv_\fU$. The desired result follows.

Combining (b),(c),(d),(e) we see that for any $y\in H_{n,D}$ we have
$$(-1)^{|\II|+|\II_\e|}\tr(\vt(y\tT_{\uD}),E^{\tca})=
\sum_{J\sub\II;\e(J)=J}(-1)^{|J_\e|}\tr(y\tT_{\uD},(\ind(E_J))^{\tca}).$$
Replacing here $(-1)^{|\II|+|\II_\e|}\vt(y\tT_{\uD})$ by $\vt(y)\tT_{\uD}$ and using Lemma 
42.6 we may rewrite this as follows:
$$\tr(\vt(y)\tT_{\uD},E^{\tca})=
\sum_{J\sub\II;\e(J)=J}(-1)^{|J_\e|}\tr(\d_J(y)\tT_{\uD},E^{\tca})$$
or equivalently (see 42.4) $\tr(\vt(y)\tT_{\uD},E^{\tca})=\tr(\d(y)\tT_{\uD},E^{\tca})$.
Since any simple $\fU(v)\ot_\ca H_n^D$-module can be obtained by extension of scalars (from
$\tca$ to $\fU(v)$) from some $E^{\tca}$ as above, we deduce that
$$\tr((\d(y)-\vt(y))\tT_{\uD},\EE)=0$$
for any simple $\fU(v)\ot_\ca H_n^D$-module $\EE$. Since $\fU(v)\ot_\ca H_n^D$ is a 
semisimple algebra, it follows that $(\d(y)-\vt(y))\tT_{\uD}$ belongs to the 
$\fU(v)$-subspace of $\fU(v)\ot_\ca H_n^D$ spanned by commutators $xx'-x'x$ with 
$x,x'\in\fU(v)\ot_\ca H_n^D$. Hence we have
$$g(\d(y)-\vt(y))\tT_{\uD}=\sum_{i=1}^mg_i
(x_i\tT_{\uD}^{s_i}x'_i\tT_{\uD}^{1-s_i}-x'_i\tT_{\uD}^{1-s_i}x_i\tT_{\uD}^{s_i})$$
with $g\in\ca-\{0\},g_i\in\ca,x_i\in H_n,x'_i\in H_n,s_i\in\ZZ$ that is,
$$g(\d(y)-vt(y))=\sum_{i=1}^mg_i
(x_i\tT_{\uD}^{s_i}x'_i\tT_{\uD}^{-s_i}-x'_i\tT_{\uD}^{1-s_i}x_i\tT_{\uD}^{s_i-1}).\tag f$$

\subhead 42.8\endsubhead
We show that for any $y,y'\in H_n$ we have

(a) $\ff\et\tio(yy'-y'\tT_{\uD}y\tT_{\uD}\i)=0$.
\nl
Let $w\in\WW,\l\in\ufs_n$. Let $\cl\in\l$. If $w\uD\l=\l$, using notation and results in 
31.6, 31.7 we have
$$\align&\ff\et\tio(v^{l(w)}\tT_w1_{\uD\l}))=gr(K^{w,\cl}_{\II,D}))\\&=
\sum_A\c^A_v(K^{w,\cl}_{\II,D}))=\sum_A\tiz^A(v^{l(w)}\tT_w1_{\uD\l}\tT_{\uD})
=\sum_A\z^A(v^{l(w)}\tT_w1_{\uD\l}\tT_{\uD})\endalign$$
(the last equation comes from 31.7(e); $A$ runs over the objects in $\hD$ up to isomorphism
such that $\z^A\ne0$.) The equation 

$\ff\et\tio(v^{l(w)}\tT_w1_{\uD\l}))=\sum_A\z^A(v^{l(w)}\tT_w1_{\uD\l}\tT_{\uD})$
\nl
holds also if $w\uD\l\ne\l$ (in this case both sides are $0$). It follows that 

$\ff\et\tio(x))=\sum_A\z^A(x\tT_{\uD})$ for any $x\in H_n$.
\nl
We deduce

$\ff\et\tio(yy'-y'\tT_{\uD}y\tT_{\uD}\i)=\sum_A(\z^A(yy'\tT_{\uD})-\z^A(y\e(y)\tT_{\uD})=0$
\nl
where the last equality follows from 31.8. This proves (a).

\proclaim{Proposition 42.9} Let $y\in H$. We have 
$\dd(\ff\et\tio(y))=\ff\et\tio(\vt(y))\in\fK(D)$ with $\dd:\fK(D)@>>>\fK(\D)$ as in 42.2.
\endproclaim
If $y\in H'_D$ (see 42.3), both sides of the desired equality are $0$. (Note that $\vt$ 
maps $H_D$ into itself and $H'_D$ into itself.) Hence we may assume that $y\in H_D$. We can
assume that $y\in H_n$ where $n\in\NN^*_\kk$. Then $y\in H_{n,D}$. By 42.4 it is enough to 
show that $\ff\et\tio(\d(y)-\vt(y))=0$. Let $g,g_i,x_i,x'_i,s_i$ be as in 
42.7(f). Since $g\ne0$ it is enough to show that $g\ff\et\tio(\d(y)-\vt(y))=0$ 
or that $\ff\et\tio(g(\d(y)-\vt(y)))=0$. Using 42.7 it is enough to show that
$$\ff\et\tio(\sum_{i=1}^mg_i
(x_i\tT_{\uD}^{s_i}x'_i\tT_{\uD}^{-s_i}-x'_i\tT_{\uD}^{1-s_i}x_i\tT_{\uD}^{s_i-1})=0.$$
Hence it is enough to show that

$\ff\et\tio(x\tT_{\uD}^sx'\tT_{\uD}^{-s}-x'\tT_{\uD}^{1-s}x\tT_{\uD}^{s-1})=0$
\nl
for any $x,x'\in H_n$ and any $s\in\ZZ$. We have
$$x\tT_{\uD}^sx'\tT_{\uD}^{-s}-x'\tT_{\uD}^{1-s}x\tT_{\uD}^{s-1}=
(z-\tT_{\uD}^{-s}z\tT_{\uD}^s)+(z'x'-x'\tT_{\uD}z'\tT_{\uD}\i)$$
where $z=x\tT_{\uD}^sx'\tT_{\uD}^{-s}\in H_n$ and $z'=\tT_{\uD}^{-s}x\tT_{\uD}^s\in H_n$. 
Hence it is enough to show that $\ff\et\tio(z'x'-x'\tT_{\uD}z'\tT_{\uD}\i)=0$ 
(see 42.8(a)) and

(a) $\ff\et\tio(z-\tT_{\uD}^{-s}z\tT_{\uD}^s)=0$
\nl
for any $z\in H_n$. This follows from 41.6(c).


\widestnumber\key{BBD}
\Refs
\ref\key{\BBD}\by A.Beilinson, J.Bernstein, P.Deligne\paper Faisceaux pervers\jour 
Ast\'erisque\vol100\yr1982\endref
\ref\key{\DE}\by P.Deligne\paper La conjecture de Weil,II\jour Publ.Math. I.H.E.S.\vol52
\yr1980\pages137-252\endref
\ref\key{\GI}\by V.Ginzburg\paper Admissible modules on a symmetric space\jour
Ast\'erisque\vol173-174\yr1989\pages199-255\endref 
\ref\key{\GR}\by I.Grojnowski\paper Character sheaves on symmetric spaces, Ph.D. thesis,
MIT\yr1992\endref
\ref\key{\CS}\by G.Lusztig\paper Character sheaves,I\jour Adv. Math.\vol56\yr1985\pages
193-237\moreref II\vol57\yr1985\pages226-265\moreref III\vol57\yr1985\pages266-315
\moreref IV\vol59\yr1986\pages1-63\moreref V\vol61\yr1986\pages103-155\endref
\ref\key{\AD}\by G.Lusztig\paper Character sheaves on disconnected groups,I\jour 
Represent. Th. (electronic)\vol7\yr2003\pages374-403\moreref II\vol8\yr2004\pages72-124
\moreref III\vol8\yr2004\pages125-144\moreref IV\vol8\yr2004\pages145-178\moreref Errata
\vol8\yr2004\pages179-179\moreref V\vol8\yr2004\pages346-376\moreref VI\vol8\yr2004\pages
377-413\moreref VII\vol9\yr2005\pages209-266\moreref VIII, math.RT/0509356\endref
\ref\key{\CRG}\by G.Lusztig\book Characters of reductive groups over a finite field,
Ann.Math.Studies\vol107\publ Princeton U.Press\yr1984\endref
\ref\key{\PCS}\by G.Lusztig\paper Parabolic character sheaves,I\jour Moscow Math.J.\vol4
\yr2004\pages153-179\endref
\ref\key{\MV}\by I.Mirkovi\'c, K.Vilonen\paper Characteristic varieties of character 
sheaves\jour Invent.Math.\vol93\yr1988\pages405-418\endref 
\endRefs
\enddocument